# UNITARY DUAL OF QUASI-SPLIT $PGSO_8^E$

CAIHUA LUO

ABSTRACT. In this paper, we first determine the explicit Langlands classification for quasi-split groups $PGSO_8^E$ by adapting Casselman–Tadić's Jacquet module machine. Based on the classification, we further sort out the unitary dual of $PGSO_8^E$ and compute the Aubert duality.

## INTRODUCTION

Let $PGSO_8^E$ be an adjoint quasi-split group of type $D_4$ over a non-archimedean field $F$ of characteristic zero, where $E$ is a cubic field extension of $F$. As part of the Langlands program, it is pivotal to understand the decomposition of induced representations and classify the unitary dual. Following Harish-Chandra, Knapp–Stein et al developed the R-group theory to determine the structure of tempered induced representations (cf. [KS71, Sil79]), and based on the R-group theory Winarsky [Win78], Keys [Key82] et al have completely determined the structure of tempered principal series for split $p$-adic Chevalley groups. As for generalized principal series (tempered or not), Shahidi [Sha90] has built up the Langlands–Shahidi theory to tackle this problem and produced quite fruitful results [Sha92, Gol94]. Along another direction, Casselman [Cas95], Rodier [Rod81], Tadić–Sally [Tad83, Tad94], Janzten [Jan96] et al have developed the Jacquet module machine to analyze the constituents of non-tempered principal series representations. But it is still far from its completeness (to my best knowledge). Motivated by the work of Rodier on regular characters, it should be reasonable to believe the existence of an internal structure for the non-tempered principal series. On the other hand, in light of unitary dual, Vogan and his collaborators have produced many influential work and created a unitary kingdom (cf. [Vog86, Vog94, KVJ16, Vog00]). As a test, some low rank groups have been computed (cf. [ST93, Mui98, Kon01, Mat10, HM10, Sch14]). From the perspective of global Langlands conjectures, AZSS (Aubert, Zelevinsky, Schneider–Stuhler) duality also plays an important role in formulating Arthur's conjecture [Art13] (as always cited as Aubert duality). It is conjectured that the AZSS duality preserves unitarity (cf. [Aub95, SS97]) and corresponds to the switch of $SL_2$-components of the A-parameter on the Galois side (see [Hir04]). In this paper, we will carry out the project for $PGSO_8^E$, whilst a similar result of the unitary dual of $Spin_8^E$ will be discussed somewhere else, and hope to finish $Sp_6$ in the near future to get a glimpse of possible internal structures of the decomposition of principal series.

Here is an outline of the paper. In the first section, we establish notation and recall some basic structure results for $PGSO_8^E$ with $E/F$ a cyclic extension and some basic representation theory facts. As the non-Galois case is almost the same, we will treat it as a remark accordingly throughout the paper. At last, we will do some basic computations for later use. In the second section, we compute the explicit Langlands classification for $PGSO_8^E$, while the last section is devoted to sorting out the unitary dual and showing the unitarizability of two isolated families $I_\alpha(1, I^\alpha(\chi_1,\chi_1^{-1})\otimes 1)$ with $\chi_{1_{|F^\times}=1}, \chi_1 \neq 1$, and $I_\beta(3, 1\otimes I^\beta(\chi_2,\chi_2^{-1}))$ with $\chi_2\circ N_{E/F}=1$, $\chi_2\neq 1$.

**Acknowledgements**. We are much indebted to Professor Wee Teck Gan for his constant help and support, and useful discussions on various topics. I would like to thank Professor Tamotsu Ikeda for discussions on the paper during a conference at IMUS, Seville, Spain.

## 1. PRELIMINARIES

Let $F$ be a non-archimedean field of characteristic zero, and $E$ be a cubic Galois field extension of $F$ with $Gal(E/F)=\langle\sigma\rangle$. Denote by $|\cdot|$ the absolute value of $F$ and by $|\cdot|\circ N_{E/F}$ the absolute value of $E$, and write $\nu_F = |\cdot| \circ det$ and $\nu_E = |\cdot| \circ N_{E/F} \circ det$. Given such an $E$, we know there is an associated adjoint quasi-split group $G = PGSO_8^E$ of type $D_4$.

Denote by $T$ a maximal torus and by $B = TU$ a Borel subgroup of $PGSO_8^E$. We know that the absolute root lattice

$$X^*(T) = \mathbb{Z}\langle \alpha_1 = e_1 - e_2, \alpha_2 = e_2 - e_3, \alpha_3 = e_3 - e_4, \alpha_4 = e_3 + e_4 \rangle,$$

and the absolute coroot lattice

$$X_*(T) = \mathbb{Z}\left\langle e_1, e_1 + e_2, \frac{1}{2}(e_1 + e_2 + e_3 - e_4), \frac{1}{2}(e_1 + e_2 + e_3 + e_4) \right\rangle.$$

Thus we may write $T$ as

$$T = H_{e_1}(t_1) H_{\frac{1}{2}(e_1+e_2+e_3-e_4)}(t_1^{\sigma}) H_{\frac{1}{2}(e_1+e_2+e_3+e_4)}(t_1^{\sigma^2}) H_{e_1+e_2}(t_2),$$

where $t_1 \in E^\times, t_2 \in F^\times$. Denote by $\alpha = \frac{1}{3}(\alpha_1 + \alpha_3 + \alpha_4)$ and $\beta = \alpha_2$. We have $\mathfrak{a}_T^* := X(T)_F \otimes_{\mathbb{Z}} \mathbb{R} = \mathbb{R}\langle \alpha, \beta \rangle$ and the positive Weyl chamber

$$C^+ = (\mathfrak{a}_T^*)_+ := \{x \in \mathfrak{a}_T^* : (x, \alpha) > 0, (x, \beta) > 0\} = \{s_1\alpha + s_2\beta : \frac{3}{2}s_2 < s_1 < 2s_2\}.$$

For any root $\gamma \in \{\alpha, \beta\}$, we denote $w_\gamma$ to be the corresponding reflection in the Weyl group $W = \langle w_\alpha, w_\beta \rangle$ of $G$. For Levi subgroups of $PGSO_8^E$, we have the following isomorphisms (see [GH06, Formula (2.28)])

$$\begin{array}{lrcl}
B = TU: & T & \xrightarrow{\sim} & E^\times \times F^\times \\
 & t & \longmapsto & (\alpha(t), \beta(t)) \\
P_\alpha = M_\alpha N_\alpha: & M_\alpha & \xrightarrow{\sim} & GL_2(E) \times F^\times / \Delta E^\times \\
 & t = (t_1, t_2) & \longmapsto & (diag(t_1, 1), t_2^{-1}) \\
P_\beta = M_\beta N_\beta: & M_\beta & \xrightarrow{\sim} & E^\times \times GL_2(F) / \Delta F^\times \\
 & t = (t_1, t_2) & \longmapsto & (t_1^{-1}, diag(t_2, 1))
\end{array}$$

Under the above realization, we have

$$w_\alpha(t_1, t_2) = (t_1^{-1}, N_{E/F}(t_1)t_2), \ w_\beta(t_1, t_2) = (t_1 t_2, t_2^{-1}).$$

And

$$T^W: \quad (t_1, t_2) \xrightarrow{w_\alpha} (t_1^{-1}, N_{E/F}(t_1)t_2) \xrightarrow{w_\beta} (t_1^{-1}N_{E/F}(t_1)t_2, N_{E/F}(t_1)^{-1}t_2^{-1}) \xrightarrow{w_\alpha} (t_1 N_{E/F}(t_1)^{-1}t_2^{-1}, N_{E/F}(t_1)t_2^2)$$

$$(t_1^{-1}, t_2^{-1}) \xleftarrow{w_\beta} (t_1^{-1}t_2^{-1}, t_2) \xleftarrow{w_\alpha} (t_1 t_2, N_{E/F}(t_1)^{-1}t_2^{-2}) \xleftarrow{w_\beta} (t_1 N_{E/F}(t_1)^{-1}t_2^{-1}, N_{E/F}(t_1)t_2^2)$$

And the Weyl group action on characters $(\chi_1, \chi_2)$ of $T \simeq E^\times \times F^\times$ is listed as follows.

$$\chi^W: \quad (\chi_1, \chi_2) \xrightarrow{w_\alpha} (\chi_1^{-1}\chi_2 \circ N_{E/F}, \chi_2), \xrightarrow{w_\beta w_\alpha} (\chi_1^{-1}\chi_1 \circ N_{E/F}\chi_2 \circ N_{E/F}^{-1}, \chi_1\chi_2^{-1}), \xrightarrow{w_\alpha w_\beta w_\alpha} (\chi_1\chi_1 \circ N_{E/F}^{-1}\chi_2 \circ N_{E/F}, \chi_1^{-1}\chi_2^2)$$

$$, (\chi_1^{-1}, \chi_1^{-1}\chi_2) \xleftarrow{w_\alpha w_\beta w_\alpha w_\beta w_\alpha}, (\chi_1\chi_2 \circ N_{E/F}^{-1}, \chi_1\chi_2^{-2}) \xleftarrow{w_\beta w_\alpha w_\beta w_\alpha}, (\chi_1\chi_1 \circ N_{E/F}^{-1}\chi_2 \circ N_{E/F}, \chi_1^{-1}\chi_2^2)$$

$$(\chi_1^{-1}, \chi_1^{-1}\chi_2) \xrightarrow{w_\beta w_\alpha w_\beta w_\alpha w_\beta w_\alpha} (\chi_1^{-1}, \chi_2^{-1})$$

For

$$(s_1, s_2) := 3s_1\alpha + s_2\beta \in X^*(T) \otimes_{\mathbb{Z}} \mathbb{C},$$

we define the associated unramified character of $T$ as

$$(t_1, t_2) \mapsto |N_{E/F}(t_1)|^{s_1}|t_2|^{s_2}.$$

For $\gamma \in \{\alpha, \beta\}$, $s_1, s_2 \in \mathbb{R}$, and $\chi_1 \times \chi_2$ unitary characters of $E^\times \times F^\times$, set

$$I^\gamma(s_1, s_2, \chi_1, \chi_2) = I^\gamma(|N_{E/F}(\cdot)|^{s_1}\chi_1 \otimes |\cdot|^{s_2}\chi_2) = Ind_T^{M_\gamma}(|N_{E/F}(\cdot)|^{s_1}\chi_1 \otimes |\cdot|^{s_2}\chi_2).$$

Similarly, we write $I(s_1, s_2, \chi_1, \chi_2) = Ind_T^G(|N_{E/F}(\cdot)|^{s_1}\chi_1 \otimes |\cdot|^{s_2}\chi_2)$ for the normalized induced representation from $B$ to $G$.

Let $R(G)$ be the Grothendieck group of admissible representations of finite length of $G$. We denote by $r_\gamma$ the normalized Jacquet functor w.r.t. $P_\gamma$, and by $r_\varnothing$ the normalized Jacquet functor w.r.t. $B$. Now we recall that (cf. [BDK86, Cas95])

$$\begin{aligned} r_\alpha(I_\alpha(\sigma)) &= \sigma + w_{3\alpha+2\beta}.\sigma + I^\alpha(w_{\alpha+\beta}.r_\varnothing(\sigma)) + I^\alpha(w_\beta.r_\varnothing(\sigma)) \\ r_\beta(I_\alpha(\sigma)) &= I^\beta(r_\varnothing(\sigma)) + I^\beta(w_{2\alpha+\beta}w_\alpha.r_\varnothing(\sigma) + I^\beta(w_{3\alpha+\beta}w_\alpha.r_\varnothing(\sigma)) \end{aligned}$$

if $\sigma \in R(M_\alpha)$, and

$$\begin{aligned} r_\beta(I_\beta(\sigma)) &= \sigma + w_{2\alpha+\beta}.\sigma + I^\beta(w_\alpha.r_\varnothing(\sigma)) + I^\beta(w_{3\alpha+\beta}.r_\varnothing(\sigma)) \\ r_\alpha(I_\beta(\sigma)) &= I^\alpha(r_\varnothing(\sigma)) + I^\alpha(w_{\alpha+\beta}w_\alpha.r_\varnothing(\sigma)) + I^\alpha(w_\beta w_\alpha.r_\varnothing(\sigma)) \end{aligned}$$

if $\sigma \in R(M_\beta)$.

We have the Aubert involution endomorphism of $R(G)$

$$D_G(\pi) = I \circ r_\varnothing(\pi) - I_\alpha \circ r_\alpha(\pi) - I_\beta \circ r_\beta(\pi) + \pi$$

It follows from [Aub95, SS97] that $\pm D_G(\pi)$ preserves irreducibility. And we have [Aub95, Theorem 1.7(3)],

$$\begin{aligned} D_G \circ I_\gamma &= I_\gamma \circ D_{M_\gamma}, \\ r_\gamma \circ D_G &= \tilde{w}_\gamma \circ D_{M_\gamma} \circ r_\gamma. \end{aligned}$$

Here $\tilde{w}_\alpha = w_{3\alpha+2\beta}$ and $\tilde{w}_\beta = w_{2\alpha+\beta}$.

Now we recall the Langlands quotient theorem and Casselman's temperedness criterion in the $PGSO_8^E$-setting (cf. [BW13, XI Proposition 2.6 and Corollary 2.7]) for later use as follows.

**Langlands quotient theorem.** *For irreducible tempered representations $\sigma$ of $GL_2$.*

*When $\chi_2$ is unitary and $s > 0$, the induced representation $Ind_{P_\alpha}^G(\nu_E^s\sigma \otimes \chi_2\nu_F^{-2s})$ has a unique irreducible quotient, i.e. the Langlands quotient $J_\alpha(s, \sigma \otimes \chi_2)$.*

*When $\chi_1$ is unitary and $s > 0$, the induced representation $Ind_{P_\beta}^G(\chi_1\nu_E^{\frac{-2s}{3}} \otimes \nu_F^s\sigma)$ has a unique irreducible quotient, i.e. the Langlands quotient $J_\beta(s, \chi_1 \otimes \sigma)$.*

*When $\chi_1, \chi_2$ are unitary and $\frac{3}{2}s_2 < 3s_1 < 2s_2$, the induced representation $I(s_1, s_2, \chi_1, \chi_2)$ has a unique irreducible quotient, i.e. the Langlands quotient $J(s_1, s_2, \chi_1, \chi_2)$.*

**Casselman's temperedness criterion.** *Suppose $\pi$ is an irreducible representation of $G$ supported on a minimal parabolic subgroup, then $\pi$ is square-integrable (resp. tempered) if and only if for any irreducible subquotient $(s_1, s_2, \chi_1, \chi_2)$ of $r_\varnothing(\pi)$ ($s_i \in \mathbb{R}$, $\chi_i$ unitary), we have*

$$(s_1, s_2) \in {}^+\mathfrak{a}_T^* = \{a\alpha + b\beta : a > 0,\ b > 0\} \ (resp.\ {}^+\bar{\mathfrak{a}}_T^*).$$

Notice that for $(s_1, s_2) = 3s_1\alpha + s_2\beta$, there exists $w \in W$ such that $(s_1, s_2)^w \in \bar{C}^+$ the closure of $C^+$, and we have $I(s_1, s_2, \chi_1, \chi_2) = I((s_1, s_2, \chi_1, \chi_2)^w)$ in $R(G)$, thus we may only need to analyze those $I(s_1, s_2, \chi_1, \chi_2)$ where $s_1$ and $s_2$ satisfy the condition that $0 \le \frac{3}{2}s_2 \le 3s_1 \le 2s_2$. To do so, we need to classify two pivotal data as follows.

**Singular character**. As the composition series of $I(s_1, s_2, \chi_1, \chi_2)$ have been determined completely by Rodier for regular characters $(s_1, s_2, \chi_1, \chi_2)$ and by Keys for unitary characters, it will be helpful to first sort out the singular characters. Recall that $\{(s_1, s_2)^w : w \in W\} =$

$$\{\pm(s_1, s_2), \pm(s_1-s_2, -s_2), \pm(2s_1-s_2, 3s_1-s_2), \pm(2s_1-s_2, 3s_1-2s_2), \pm(s_1-s_2, 3s_1-2s_2), \pm(s_1, 3s_1-s_2)\}.$$

So for those $s_1$ and $s_2$ satisfying the condition that $0 \le \frac{3}{2}s_2 \le 3s_1 \le 2s_2$, the set $S$ of singular characters consists of those unitary $\chi$ with multiplicity $m > 2$:

$(1, 1; D_6)$, $(\chi_1, 1; S_3; \chi_1 \ne 1,\ \chi_1|_{F^\times} = 1)$, $(\chi_1, 1; \langle w_\alpha, w_{3\alpha+2\beta}\rangle\,; \chi_1^2 = 1,\ \chi_1|_{F^\times} \ne 1)$, $(1, \chi_2; \langle w_\beta, w_{2\alpha+\beta}\rangle\,; \chi_2^2 = 1,\ \chi_2 \ne 1)$, $(\chi_1, \chi_1; \langle w_{3\alpha+\beta}, w_{\alpha+\beta}\rangle\,; \chi_1^2 = 1,\ \chi_1|_{F^\times} \ne 1)$,

and those unitary $\chi$ with multiplicity $m = 2$:

$(\chi_1,\chi_2;\langle w_\alpha w_{3\alpha+2\beta}\rangle;\chi_1^2=1,\ \chi_2^2=1)$, $(\chi_1,\chi_2;\langle w_\alpha\rangle;\chi_1^2=\chi_2\circ N_{E/F})$, $(\chi_1,\chi_2;\langle w_{\alpha+\beta}\rangle;\chi_1=\chi_2\circ N_{E/F})$, $(1,\chi_2;\langle w_{2\alpha+\beta}\rangle)$, $(\chi_1,1;\langle w_{3\alpha+2\beta}\rangle)$, $(\chi_1,\chi_2;\langle w_\beta\rangle;\chi_1=\chi_2^2)$, $(\chi_1,\chi_2;\langle w_{3\alpha+\beta}\rangle;\chi_1=\chi_2)$,
and those non-unitary $\chi$ with $m=2$:

$$(s_1,2s_1,\chi_1,\chi_2;\langle w_\alpha\rangle;s_1>0,\ \chi_1^2=\chi_2\circ N_{E/F}),\quad (s_1,\frac{3}{2}s_1,\chi_1,\chi_2;\langle w_\beta\rangle;s_1>0,\ \chi_1|_{F^\times}=\chi_2^2).$$

**Reducibility point.** In what follows, we will describe the rank 1 reducibility points, as we believe in most cases, rank 1 irreducibility should determine the irreducibility of the full induced representation. The set $R$ of rank 1 reducibility points is listed as follows.

| | Rank 1 reducibility $(s_1,s_2,\chi_1,\chi_2)$ |
|---|---|
| $\alpha^\vee$ | $2s_1-s_2=1,\ \chi_1^2=\chi_2\circ N_{E/F}$ |
| $\beta^\vee$ | $-3s_1+2s_2=1,\ \chi_1=\chi_2^2$ |
| $(\alpha+\beta)^\vee$ | $-s_1+s_2=1,\ \chi_1=\chi_2\circ N_{E/F}$ |
| $(2\alpha+\beta)^\vee$ | $s_1=1,\ \chi_1=1$ |
| $(3\alpha+\beta)^\vee$ | $3s_1-s_2=1,\ \chi_1=\chi_2$ |
| $(3\alpha+2\beta)^\vee$ | $s_2=1,\ \chi_2=1$ |

Table 1. Reducibility point

| | Reducible coroot $(\cdot)^\vee$ | relation for $(\cdot)^\vee$ | $(s_1,s_2,\chi_1,\chi_2)$ |
|---|---|---|---|
| $\#R=4$ | $(\alpha+\beta),\ \beta,\ (2\alpha+\beta),\ (3\alpha+\beta)$ | $w_\alpha$ | (1,2,1,1) |
| $\#R=2$ | $(3\alpha+\beta),\ (3\alpha+2\beta)$ | $w_\beta$ | $(\frac{2}{3},1,\chi_1,1;\chi_1|_{F^\times}=1)$ |
| | $\alpha,\ (\alpha+\beta)$ | $w_\beta$ | $(2,3,1,\chi_2;\chi_2\circ N_{E/F}=1)$ |
| | $\alpha,\ \beta$ | | $(3,5,1,1)$ |
| | $(\alpha+\beta),\ (2\alpha+\beta)$ | $w_\alpha$ | $(1,2,1,\chi_2;\chi_2\neq1\ \&\ \chi_2\circ N_{E/F}=1)$ |
| | $\beta,\ (2\alpha+\beta)$ | | $(1,2,1,\chi_2;\chi_2\neq1\ \&\ \chi_2^2=1)$ |
| | $(\alpha+\beta),\ (3\alpha+\beta)$ | | $(1,2,\chi_1,\chi_2;\chi_2\neq1\ \&\ \chi_1^2=1,\ \chi_1=\chi_2)$ |
| | $\beta,\ (3\alpha+\beta)$ | $w_\alpha$ | $(1,2,\chi_1,1;\chi_1\neq1\ \&\ \chi_1|_{F^\times}=1)$ |

Table 2. $\#R>1$

Before moving to the next computation section, we recall Shahidi's local coefficient formula in the $PGSO_8^E$-setting based on its multiplicative property as follows (please refer to [Sha90] for the notions), up to a monomial in $q^{-s}$,

$$\begin{aligned}
C_\psi(s,\delta(\chi_1)\otimes\chi_2,w_{3\alpha+2\beta}) &= \frac{L_F(\frac{5}{2}-s,\chi_2^2\chi_1)L_F(1-2s,\chi_2)L_F(-\frac{1}{2}-s,\chi_1^{-1}\chi_2^{-1})}{L_F(s-\frac{3}{2},\chi_2^{-2}\chi_1^{-1})L_F(2s,\chi_2^{-1})L_F(\frac{3}{2}+s,\chi_1\chi_2)}\\
&\times\frac{L_E(\frac{3}{2}-s,\chi_1^{-2}(\chi_1\chi_2)\circ N_{E/F})L_E(\frac{1}{2}-s,\chi_1^{-1}\circ N_{E/F}\chi_1^2)}{L_E(s-\frac{1}{2},\chi_1^2(\chi_1\chi_2)^{-1}\circ N_{E/F})L_E(s+\frac{1}{2},\chi_1\circ N_{E/F}\chi_1^{-2})}\\
C_\psi(s,\chi_1\otimes\delta(\chi_2),w_{2\alpha+\beta}) &= \frac{L_E(\frac{3}{2}-\frac{s}{3},\chi_1^2\chi_2\circ N_{E/F})L_E(1-\frac{2s}{3},\chi_2^{-2}\chi_1\circ N_{E/F})L_E(\frac{1}{2}-\frac{s}{3},\chi_1^2(\chi_1\chi_2)^{-1}\circ N_{E/F})}{L_E(-\frac{1}{2}+\frac{s}{3},\chi_1^{-2}\chi_2^{-1}\circ N_{E/F})L_E(\frac{2s}{3},\chi_2\chi_1^{-1}\circ N_{E/F})L_E(\frac{1}{2}+\frac{s}{3},\chi_1^{-2}(\chi_1\chi_2)\circ N_{E/F})}\\
&\times\frac{L_F(\frac{3}{2}-s,\chi_1\chi_2)L_F(\frac{1}{2}-s,\chi_2^{-1})}{L_F(s-\frac{1}{2},(\chi_1\chi_2)^{-1})L_F(\frac{1}{2}+s,\chi_2)}
\end{aligned}$$

In view of the above formulas, we have the following lemma which results from [Sha81, Proposition 3.3.1].

**Lemma 1.1.** *We have the genericity of those representations which will be used in the next section.*

$$J_\alpha(\nu_E^{\frac{1}{2}}\delta(1)\otimes\chi_2\nu_F^{-1})|_{\chi_2\neq1,\chi_2\circ N_{E/F}=1},\quad J_\alpha(\nu_E^{\frac{3}{2}}\delta(1)\otimes\nu_F^{-3}),\quad J_\beta(\nu_E^{\frac{-1}{3}}\chi_1\otimes\nu_F^{\frac{1}{2}}\delta(1))|_{\chi_1|_{F^\times}=1}.$$

## 2. Langlands classification

In this section, we will carry out the computation of the constituents of principal series in detail following Casselman–Tadić's Jacquet module machine. Recall that given a character $\chi := (s_1, s_2, \chi_1, \chi_2)$ of $T$, under the previous realization of Levi subgroups, we have

$$I^\alpha(\chi) = Ind^{GL_2}(|N_{E/F}(\cdot)|^{s_1}\chi_1 \otimes |N_{E/F}(\cdot)|^{-s_1+s_2}\chi_1^{-1}\chi_2 \circ N_{E/F}) \otimes \chi_2^{-1}|\cdot|^{-s_2},$$
$$I^\beta(\chi) = \chi_1^{-1}|N_{E/F}(\cdot)|^{-s_1} \otimes Ind^{GL_2}(|\cdot|^{s_2}\chi_2 \otimes |\cdot|^{3s_1-s_2}\chi_1\chi_2^{-1}).$$

It is well-known that they are reducible if and only if

$$(2s_1 - s_2, \chi_1^2\chi_2^{-1} \circ N_{E/F}) = (\pm 1, 1) \text{ and } (2s_2 - 3s_1, \chi_2^2\chi_1^{-1}) = (\pm 1, 1) \text{ respectively.}$$

And their Jacquet modules $r_\varnothing$ have the form

$$\begin{aligned}
r_\varnothing^{M_\alpha}(I^\alpha(\chi)) =& \{(s_1, s_2, \chi_1, \chi_2), (-s_1 + s_2, s_2, \chi_1^{-1}\chi_2 \circ N_{E/F}, \chi_2)\},\\
r_\varnothing^{M_\beta}(I^\beta(\chi)) =& \{(s_1, s_2, \chi_1, \chi_2), (s_1, 3s_1 - s_2, \chi_1, \chi_1\chi_2^{-1})\},\\
r_\varnothing^{G}(I^G(\chi)) \overset{M_\alpha}{=}& \{\pm(s_1, s_2, \chi_1, \chi_2), \pm(-s_1 + s_2, s_2, \chi_1^{-1}\chi_2 \circ N_{E/F}, \chi_2)\}\\
&\cup\{\pm(s_1, 3s_1 - s_2, \chi_1, \chi_1\chi_2^{-1}), \pm(2s_1 - s_2, 3s_1 - s_2, \chi_1^{-1}\chi_1 \circ N_{E/F}\chi_2^{-1} \circ N_{E/F}, \chi_1\chi_2^{-1})\}\\
&\cup\{\pm(s_1 - s_2, 3s_1 - 2s_2, \chi_1\chi_2^{-1} \circ N_{E/F}, \chi_1\chi_2^{-2}),\\
&\quad \pm(2s_1 - s_2, 3s_1 - 2s_2, \chi_1^{-1}\chi_1 \circ N_{E/F}\chi_2^{-1} \circ N_{E/F}, \chi_1\chi_2^{-2})\},\\
r_\varnothing^{G}(I^G(\chi)) \overset{M_\beta}{=}& \{\pm(s_1, s_2, \chi_1, \chi_2), \pm(s_1, 3s_1 - s_2, \chi_1, \chi_1\chi_2^{-1})\}\\
&\cup\{\pm(s_1 - s_2, 3s_1 - 2s_2, \chi_1\chi_2^{-1} \circ N_{E/F}, \chi_1\chi_2^{-2}), \mp(-s_1 + s_2, s_2, \chi_1^{-1}\chi_2 \circ N_{E/F}, \chi_2)\}\\
&\cup\{\pm(2s_1 - s_2, 3s_1 - 2s_2, \chi_1^{-1}\chi_1 \circ N_{E/F}\chi_2^{-1} \circ N_{E/F}, \chi_1\chi_2^{-2}),\\
&\quad \pm(2s_1 - s_2, 3s_1 - s_2, \chi_1^{-1}\chi_1 \circ N_{E/F}\chi_2^{-1} \circ N_{E/F}, \chi_1\chi_2^{-1})\}.
\end{aligned}$$

Now suppose that $s_1$ and $s_2$ satisfy the condition that $2s_2 \geq 3s_1 \geq \frac{3}{2}s_2 \geq 0$, we are ready to carry out the tedious but excited computation case by case as follows, as it may show some hidden structures.

<u>$\#R = 0$, $(s_1, s_2, \chi_1, \chi_2)$ non-unitary</u>. **Claim**: $I(\chi)$ is irreducible.

If $\chi$ is regular, i.e. $Stab_W(\chi) = 1$, the diagram chasing looks pretty easy. We write down the diagram as a template for other cases.

$$\begin{aligned}
&\overset{\alpha}{\Longrightarrow} \{(s_1, s_2, \chi_1, \chi_2), (-s_1 + s_2, s_2, \chi_1^{-1}\chi_2 \circ N_{E/F}, \chi_2)\}\\
&\overset{\beta}{\Longrightarrow} \{(s_1, s_2, \chi_1, \chi_2), (s_1, 3s_1 - s_2, \chi_1, \chi_1\chi_2^{-1})\};\\
&\qquad \{-(s_1 - s_2, 3s_1 - 2s_2, \chi_1\chi_2^{-1} \circ N_{E/F}, \chi_1\chi_2^{-2}), (-s_1 + s_2, s_2, \chi_1^{-1}\chi_2 \circ N_{E/F}, \chi_2)\}\\
&\overset{\alpha}{\Longrightarrow} \{(s_1, 3s_1 - s_2, \chi_1, \chi_1\chi_2^{-1}), (2s_1 - s_2, 3s_1 - s_2, \chi_1^{-1}\chi_1 \circ N_{E/F}\chi_2^{-1} \circ N_{E/F}, \chi_1\chi_2^{-1})\};\\
&\qquad \{-(s_1 - s_2, 3s_1 - 2s_2, \chi_1\chi_2^{-1} \circ N_{E/F}, \chi_1\chi_2^{-2}), -(2s_1 - s_2, 3s_1 - 2s_2, \chi_1^{-1}\chi_1 \circ N_{E/F}\chi_2^{-1} \circ N_{E/F}, \chi_1\chi_2^{-2})\}\\
&\overset{\beta}{\Longrightarrow} \{(2s_1 - s_2, 3s_1 - 2s_2, \chi_1^{-1}\chi_1 \circ N_{E/F}\chi_2^{-1} \circ N_{E/F}, \chi_1\chi_2^{-2}), (2s_1 - s_2, 3s_1 - s_2, \chi_1^{-1}\chi_1 \circ N_{E/F}\chi_2^{-1} \circ N_{E/F}, \chi_1\chi_2^{-1})\};\\
&\qquad \{-(2s_1 - s_2, 3s_1 - 2s_2, \chi_1^{-1}\chi_1 \circ N_{E/F}\chi_2^{-1} \circ N_{E/F}, \chi_1\chi_2^{-2}), -(2s_1 - s_2, 3s_1 - s_2, \chi_1^{-1}\chi_1 \circ N_{E/F}\chi_2^{-1} \circ N_{E/F}, \chi_1\chi_2^{-1})\}\\
&\overset{\alpha}{\Longrightarrow} \{(s_1 - s_2, 3s_1 - 2s_2, \chi_1\chi_2^{-1} \circ N_{E/F}, \chi_1\chi_2^{-2}), (2s_1 - s_2, 3s_1 - 2s_2, \chi_1^{-1}\chi_1 \circ N_{E/F}\chi_2^{-1} \circ N_{E/F}, \chi_1\chi_2^{-2})\};\\
&\qquad \{-(s_1, 3s_1 - s_2, \chi_1, \chi_1\chi_2^{-1}), -(2s_1 - s_2, 3s_1 - s_2, \chi_1^{-1}\chi_1 \circ N_{E/F}\chi_2^{-1} \circ N_{E/F}, \chi_1\chi_2^{-1})\}\\
&\overset{\beta}{\Longrightarrow} \{(s_1 - s_2, 3s_1 - 2s_2, \chi_1\chi_2^{-1} \circ N_{E/F}, \chi_1\chi_2^{-2}), -(-s_1 + s_2, s_2, \chi_1^{-1}\chi_2 \circ N_{E/F}, \chi_2)\};\\
&\qquad \{-(s_1, s_2, \chi_1, \chi_2), -(s_1, 3s_1 - s_2, \chi_1, \chi_1\chi_2^{-1})\}\\
&\overset{\alpha}{\Longrightarrow} \{-(s_1, s_2, \chi_1, \chi_2), -(-s_1 + s_2, s_2, \chi_1^{-1}\chi_2 \circ N_{E/F}, \chi_2)\}.
\end{aligned}$$

Whence $I(\chi)$ is irreducible. If $\chi$ is singular, as the singularity is given by $\langle w_\alpha \rangle$ or $\langle w_\beta \rangle$, we may obtain $I(\chi)$ is irreducible as well by the same argument.

<u>$(1, 2, 1, 1; \langle w_\alpha \rangle)$, $(\#R = 4, w_\alpha)$</u>. **Claim**: $I(\chi)$ is of length $2^{\#R/2} + 2$ and multiplicity at most 2, and the two subrepresentations are square-integrable.

Comparing Table 1 and Table 2, we find that $\chi$ is singular. The Jacquet modules $r_\varnothing$ are listed as follows, write $(s_1, s_2)$ for $(s_1, s_2, 1, 1)$ for simplicity.

The subrepresentation $\pi(1)$:

$$2(1,2),\ (1,1);$$

Subrepresentation $\pi(1)'$:

$$(1,1);$$

Subquotient $J_\alpha(\nu_E^{1/2}\delta(1)\otimes\nu_F^{-1})$ (multiplicity 2):

$$(0,1),\ (0,-1);$$

Quotient $J_\beta(\nu_E^{-1}\otimes\nu_F^{3/2}\delta(1))$:

$$(-1,-1);$$

The Langlands quotient $J_\alpha(\nu_E I^\alpha(1\otimes 1)\otimes\nu_F^{-2})$:

$$2(-1,-2),\ (-1,-1).$$

*Proof.* In $R(G)$,

$$I(1,2)=I(0,1)=I(1,1)=I_\alpha(\nu_E^{1/2}\delta(1)\otimes\nu_F^{-1})+I_\alpha(\nu_E^{1/2}1_{GL_2}\otimes\nu_F^{-1})=I_\beta(\nu_E^{-1}\otimes\nu_F^{3/2}1_{GL_2})+I_\beta(\nu_E^{-1}\otimes\nu_F^{3/2}\delta(1)).$$

We write the semisimplification of Jacquet modules as follows.

$$r_\beta(I_\alpha(\nu_E^{1/2}\delta(1)\otimes\nu_F^{-1}))=2\{(1,1)\}+2\{(1,2)\}+\{(0,1),(0,-1)\}$$
$$r_\beta(I_\beta(\nu_E^{-1}\otimes\nu_F^{3/2}1_{GL_2}))=\{(1,1)\}+2\{(-1,-2)\}+\{(-1,-1)\}+\{(0,1),(0,-1)\}$$

It is easy to see

$$\pi(1)':=I_\alpha(\nu_E^{1/2}\delta(1)\otimes\nu_F^{-1})\cap I_\beta(\nu_E^{-1}\otimes\nu_F^{3/2}1_{GL_2})\neq\varnothing.$$

Notice that

$$J_\alpha(\nu_E^{1/2}\delta(1)\otimes\nu_F^{-1})\hookrightarrow I_\alpha(\nu_E^{-1/2}\delta(1)\otimes\nu_F)\hookrightarrow I(0,-1),$$

it implies $r_\varnothing(\pi(1)')=(1,1)$.

Now consider

$$I(0,-1)\simeq I(0,1)=I_\alpha(\nu_E^{1/2}1_{GL_2}\otimes\nu_F^{-1})+I_\alpha(\nu_E^{1/2}\delta(1)\otimes\nu_F^{-1}).$$

We write the semisimplification of Jacquet modules as follows.

$$r_\beta(I_\alpha(\nu_E^{-1/2}\delta(1)\otimes\nu_F))=2\{(1,1)\}+2\{(1,2)\}+\{(0,1),(0,-1)\}$$
$$r_\beta(I_\alpha(\nu_E^{1/2}1_{GL_2}\otimes\nu_F^{-1}))=2\{(-1,-2)\}+2\{(-1,-1)\}+\{(0,1),(0,-1)\}$$

It is easy to see

$$I_\alpha(\nu_E^{-1/2}\delta(1)\otimes\nu_F)\cap I_\alpha(\nu_E^{1/2}1_{GL_2}\otimes\nu_F^{-1})\neq\varnothing$$

with the Jacquet module $\{(0,1),(0,-1)\}$.

Note also that under the Aubert duality,

$$r_\varnothing\circ D_G(\pi(1)')=(-1,-1).$$

Observe that

$$r_\beta\circ I_\beta(\nu_E^{-1}\otimes\nu_F^{3/2}\delta(1))=2\{(1,2)\}+\{(0,1),(0,-1)\}+\{(1,1)\}+\{(-1,-1)\}$$

and the possible Langlands quotients associated to $I(1,2)$ are

$$J_\beta(\nu_E^{-1}\otimes\nu_F^{3/2}\delta(1)),\ J_\alpha(\nu_E^{1/2}\delta(1)\otimes\nu_F^{-1})\ \textit{and}\ J_\alpha(\nu_E I^\alpha(1\otimes 1)\otimes\nu_F^{-2}),$$

we may conclude that $D_G(\pi(1)')$ is of multiplicity one in $I(1,2)$. □

$\underline{(1,2,1,\chi_2;\langle w_\alpha\rangle;\chi_2\neq 1 \ \&\ \chi_2\circ N_{E/F}=1),\ (\#R=2,w_\alpha)}$. **Claim**: $I(\chi)$ is of length $2^{\#R/2}$ and multiplicity 1, and the subrepresentation is square-integrable and maps to the Langlands quotient under the Aubert duality.

Comparing Table 1 and Table 2, we find that $\chi$ is singular. The above claim follows from the fact that $J_\alpha(\nu_E^{1/2}\delta(1)\otimes\chi_2\nu_F^{-1})$ is generic and Rodier's heredity theorem [Rod73, Theorem 2]. The Jacquet modules $r_\varnothing$ are listed as follows.

The subrepresentation $I_\alpha(\nu_E^{1/2}\delta(1)\otimes\chi_2\nu_F^{-1})$:

$$2(1,2,1,\chi_2),2(1,1,1,\chi_2^{-1}),(0,1,1,\chi_2^2),(0,-1,1,\chi_2);$$

The Langlands quotient $J_\alpha(\nu_E I^\alpha(1\otimes 1)\otimes\chi_2^{-1}\nu_F^{-2})$:

$$-2(1,2,1,\chi_2),-2(1,1,1,\chi_2^{-1}),(0,1,1,\chi_2^2),(0,-1,1,\chi_2).$$

$\underline{(1,2,1,\chi_2;\langle w_\alpha\rangle;\chi_2\neq 1 \ \&\ \chi_2^2=1),\ (\#R=2,\varnothing)}$. **Claim**: $I(\chi)$ is of length $2^{\#R}$ and multiplicity 1, and the subrepresentation is square-integrable and maps to the Langlands quotient under the Aubert duality.

Comparing Table 1 and Table 2, we find that $\chi$ is regular. The Jacquet modules $r_\varnothing$ are listed as follows.

The subrepresentation $\pi(\chi_2)$:

$$(1,2,1,\chi_2),\ (1,2,\chi_2\circ N_{E/F},\chi_2),\ (1,1,\chi_2\circ N_{E/F},1);$$

Subquotient $J_\beta(\nu_E^{-1}\otimes\nu_F^{3/2}\delta(\chi_2))$:

$$(0,-1,\chi_2\circ N_{E/F},\chi_2^2),\ (0,1,\chi_2\circ N_{E/F},\chi_2),\ (-1,-1,1,\chi_2);$$

Subquotient $J_\alpha(\nu_E^{1/2}\delta(\chi_2\circ N_{E/F})\otimes\nu_F^{-1})$:

$$-(0,-1,\chi_2\circ N_{E/F},\chi_2^2),\ -(0,1,\chi_2\circ N_{E/F},\chi_2),\ -(-1,-1,1,\chi_2);$$

The Langlands quotient $J_\alpha(\nu_E I^\alpha(1\otimes\chi_2\circ N_{E/F})\otimes\chi_2^{-1}\nu_F^{-2})$:

$$-(1,2,1,\chi_2),\ -(1,2,\chi_2\circ N_{E/F},\chi_2),\ -(1,1,\chi_2\circ N_{E/F},1).$$

$\underline{(1,2,\chi_1,\chi_2;\langle w_\alpha\rangle;\chi_2\neq 1 \ \&\ \chi_1^2=1,\ \chi_1=\chi_2),\ (\#R=2,\varnothing)}$. **Claim**: $I(\chi)\simeq I(1,2,1,\chi_2;\chi_2^2=1)$.

As $I(1,2,1,\chi_2)=$

$$I_\alpha(I^\alpha(\nu_E\otimes\nu_E\chi_2\circ N_{E/F})\otimes\chi_2\nu_F^{-2})\simeq I_\alpha(I^\alpha(\nu_E\chi_2\circ N_{E/F}\otimes\nu_E)\otimes\chi_2\nu_F^{-2})=I(1,2,\chi_2\circ N_{E/F},\chi_2).$$

$\underline{(1,2,\chi_1,1;\chi_1\neq 1 \ \&\ \chi_1|_{F^\times}=1),\ (\#R=2,w_\alpha)}$. **Claim**: $I(\chi)\simeq I(1,2,\chi_1^{-1},1)$ is of length $2^{\#R}$ and multiplicity 1, and the subrepresentation is square-integrable and maps to the Langlands quotient under the Aubert duality.

Comparing Table 1 and Table 2, we find that $\chi$ is regular. The Jacquet modules $r_\varnothing$ are listed as follows, write $(s_1,s_2,\mu)$ for $(s_1,s_2,\mu,1)$ for simplicity.

The subrepresentation $\pi(\chi)$:

$$(1,2,\chi_1),\ (1,2,\chi_1^{-1});$$

Subquotient $J_\beta(\nu_E^{-1}\chi_1\otimes\nu_F^{3/2}\delta(1))$:

$$(1,1,\chi_1),\ (0,1,\chi_1^{-1}),\ (0,-1,\chi_1^{-1}),\ (-1,-1,\chi_1);$$

Subquotient $J_\beta(\nu_E^{-1}\chi_1^{-1}\otimes\nu_F^{3/2}\delta(1))$:

$$-(1,1,\chi_1),\ -(0,1,\chi_1^{-1}),\ -(0,-1,\chi_1^{-1}),\ -(-1,-1,\chi_1);$$

The Langlands quotient $J_\alpha(\nu_E I^\alpha(\chi_1,\chi_1^{-1})\otimes\nu_F^{-2})$:

$$(-1,-2,\chi_1^{-1}),\ (-1,-2,\chi_1).$$

$\underline{(2/3,1,\chi_1,1;\langle w_\beta\rangle;\chi_1|_{F^\times}=1),\ (\#R=2,w_\beta)}$. **Claim**: $I(\chi)$ is of length $2^{\#R/2}$ and multiplicity 1, and the subrepresentation maps to the Langlands quotient under the Aubert duality.

The above claim follows from the fact that $J_\beta(\nu_E^{-1/3}\chi_1\otimes\nu_F^{1/2}\delta(1))$ is generic and Rodier's heredity theorem [Rod73, Theorem 2]. The Jacquet modules $r_\varnothing$ are listed as follows, write $(s_1,s_2)$ for $(s_1,s_2,1,1)$ for simplicity.

The subrepresentation $I_\beta(\nu_E^{-1/3}\chi_1\otimes\nu_F^{1/2}\delta(1))$:

$$2(2/3,1),\ 2(1/3,1),\ (1/3,0),\ (-1/3,0);$$

The Langlands quotient $J_\beta(\nu_E^{-2/3}\chi_1^{-1}\otimes\nu_F I^\beta(1\otimes1))$:

$$2(-2/3,-1),\ 2(-1/3,-1),\ (1/3,0),\ (-1/3,0).$$

$\underline{(2,3,1,\chi_2;\langle w_\beta\rangle;\chi_2\circ N_{E/F}=1,\chi_2\neq1),\ (\#R=2,w_\beta)}$. **Claim**: $I(\chi)\simeq I(2,3,1,\chi_2^{-1})$ is of length $2^{\#R}$ and multiplicity 1, and the subrepresentation is square-integrable and maps to the Langlands quotient under the Aubert duality.

It is easy to see that such a $\chi$ is regular. The Jacquet modules are listed as follows.

The subrepresentation $\pi(\chi_2)$:

$$(2,3,1,\chi_2),\ (2,3,1,\chi_2^{-1});$$

Subquotient $J_\alpha(\nu_E^{3/2}\delta(1)\otimes\chi_2\nu_F^{-3})$:

$$(1,3,1,\chi_2),\ (1,0,1,\chi_2^{-1}),\ (-1,0,1,\chi_2^{-1}),\ (-1,-3,1,\chi_2);$$

Subquotient $J_\alpha(\nu_E^{3/2}\delta(1)\otimes\chi_2^{-1}\nu_F^{-3})$:

$$-(1,3,1,\chi_2),\ -(1,0,1,\chi_2^{-1}),\ -(-1,0,1,\chi_2^{-1}),\ -(-1,-3,1,\chi_2);$$

The Langlands quotient $J_\beta(\nu_E^{-2}\otimes\nu_F^{3}I^\beta(\chi_2\otimes\chi_2^{-1}))$:

$$-(2,3,1,\chi_2),\ -(2,3,1,\chi_2^{-1}).$$

$\underline{(2,3,1,1;\langle w_\beta\rangle),\ (\#R=2,w_\beta)}$. **Claim**: $I(\chi)$ is of length $2^{\#R/2}$ and multiplicity 1, and the subrepresentation maps to the Langlands quotient under the Aubert duality.

It is easy to see that such a $\chi$ is singular. The above claim follows from the fact that $J_\alpha(\nu_E^{3/2}\delta(1)\otimes\nu_F^{-3})$ is generic and Rodier's heredity theorem [Rod73, Theorem 2]. The Jacquet modules $r_\varnothing$ are listed as follows, write $(s_1,s_2)$ for $(s_1,s_2,1,1)$ for simplicity.

The subrepresentation $I_\alpha(\nu_E^{3/2}\delta(1)\otimes\chi_2^{-1}\nu_F^{-3})$:

$$2(2,3),\ (1,0),\ (1,3),\ (-1,0),\ (-1,-3);$$

The Langlands quotient $J_\beta(\nu_E^{-2}\otimes\nu_F^{3}I^\beta(\chi_2\otimes\chi_2^{-1}))$:

$$2(-2,-3),\ (1,0),\ (1,3),\ (-1,0),\ (-1,-3).$$

$\underline{(3,5,1,1),(\#R=2)}$. **Claim**: $I(\chi)$ is of length $2^{\#R}$ and multiplicity 1, and the subrepresentation is square-integrable and maps to the Langlands quotient under the Aubert duality.

The Jacquet modules $r_\varnothing$ of the constituents of $I(\chi)$ are listed as follows, write $(s_1,s_2)$ for $(s_1,s_2,1,1)$ for simplicity.

The subrepresentation $St_G$:

$$(3,5);$$

Subquotient $J_\beta(\nu_E^{-3}\otimes\nu_F^{9/2}\delta(1))$:

$$(2,5),\ (2,1),\ (1,-1),\ (-1,-4),\ (-3,-4);$$

Subquotient $J_\alpha(\nu_E^{5/2}\delta(1)\otimes\nu_F^{-5})$:

$$(3,4),\ (1,4),\ (1,-1),\ (-2,-1),\ (-2,-5);$$

The Langlands quotient $1_G$:

$$(-3,-5).$$

$\#R=1$. **Claim**: $I(\chi)$ is of length 2 and multiplicity 1, and the subrepresentation maps to the quotient under the Aubert duality.

Comparing Table 1 and Table 2, we find that only

$$(1,3/2,1,\chi_2;\langle w_\beta\rangle;\chi_2^2=1) \text{ and } (1/2,1,\chi_1,1;\langle w_\alpha\rangle;\chi_1^2=1)$$

are singular characters. The claim that $I(\chi)$ is of length 2 can be checked easily by diagram chasing. As for multiplicity 1, notice that the singularity given by $\langle w_\alpha\rangle$ or $\langle w_\beta\rangle$ is not the one giving rise to the rank 1 reducibility, so it is of multiplicity one.

$\#R=0,\ (\chi_1,\chi_2;\langle ?\rangle)|_{?^2=1}$. **Claim**: $I(\chi)$ is irreducible except the case $?=w_\alpha w_{3\alpha+2\beta}$ which is reducible (see the paragraph below) and its constituents are invariant under the Aubert duality.

For $\langle w_\alpha w_{3\alpha+2\beta}\rangle$, if $I(\chi)$ is reducible, then it is of multiplicity one, otherwise $dim\, Hom(I(\chi),I(\chi))\le 2$, contradiction.

As for other cases, it is easy to verify that $I(\chi)$ is irreducible.

Other $(\chi_1,\chi_2)$. **Claim**: They are irreducible.

Note that the rank 1 groups are $GL_2(F)\times E^\times/\Delta F^\times$ and $GL_2(E)\times F^\times/\Delta E^\times$, so the associated Plancherel measures of unitary induced representations are the same as in $GL_2(F)$ and $GL_2(E)$ respectively. So by Keys' theorem [Key82, Theorem 1], the $R$-group can be described as

$$R=\{w\in W_\chi=Stab_W(\chi):\ \gamma>0 \text{ and } \chi_\gamma:=\chi\circ\gamma^\vee=1 \text{ imply that } w.\gamma>0\}.$$

Whence they are reducible unless $\chi_1$ and $\chi_2$ are different characters of order 2 which results from the same reason as in [Key82, Theorem $G_2$].

*Remark* **1.** *From the above computation for the case $(\#R=2,m=2)$, we know that $I(\chi)$ is of length 2 and multiplicity one. Heuristically, this may be a general result for reductive groups based on the following strategy by a case-by-case check.*

*(i) Possible Jacquet module decomposition of $r_\varnothing(I(\chi))$:*

$$\{\chi^w:\ w\in W_1Stab_W(\chi)\},\quad \{\chi^w:\ w\in W_2\},\quad \{\chi^w:\ w\in W_2\},\quad \{\chi^w:\ w\in W_3Stab_W(\chi)\},$$

*where $W_i,\ i=1,2,3,$ are subsets of the Weyl group $W$.*

*(ii) Genericity of the quotient $\pi$ of the subrepresentation of $I(\chi)$ associated to $\{\chi^w:\ w\in W_2\}$. This may be checked using the Langlands–Shahidi theory.*

*(iii) Rodier's heredity theorem which implies that the generic subquotient of $I(\chi)$ is of multiplicity one.*

*Note that once we know $\pi$ is generic, the above assertion also follows from the standard module conjecture proved by Muić and Heiermann.*

**Corollary 2.1.** *(i) $I_\alpha(s,\delta(\chi_1)\otimes\chi_2)$ reduces if and only if*

$$s=\pm1/2,\chi_1=\chi_1\circ N_{E/F},\chi_2=1 \quad or \quad s=\pm3/2,\chi_1=1,\chi_2\neq1,\chi_2\circ N_{E/F}=1 \quad or \quad s=\pm5/2,\chi_1=1,\chi_2=1;$$

*(ii) $I_\beta(s,\chi_1\otimes\delta(\chi_2))$ reduces if and only if*

$$s=\pm3/2,\chi_1|_{F^\times}=1,\chi_2=1 \quad or \quad s=\pm3/2,\chi_1=1,\chi_2^2=1 \quad or \quad s=\pm9/2,\chi_1=1,\chi_2=1.$$

**Conclusion**. In what follows, we summarize our previous computation for later use.

| Regular $\#R=1$ | $I(s_1,s_2,\chi_1,\chi_2)$ | |
|---|---|---|
| | subrepresentation | Langlands quotient |
| $2s_1-s_2=1$ and $\chi_1^2=\chi_2\circ N_{E/F}$ | $I_\alpha(s_1-1/2,\delta(\chi_1)\otimes\chi_2^{-1})$ | $I_\alpha(s_1-1/2,\chi_1\circ det\otimes\chi_2^{-1})$ |
| $2s_2-3s_1=1$ and $\chi_2^2=\chi_1$ | $I_\beta(s_2-1/2,\chi_1^{-1}\otimes\delta(\chi_2))$ | $I_\beta(s_2-1/2,\chi_1^{-1}\otimes\chi_2\circ det)$ |
| $s_2-s_1=1$ and $\chi_1=\chi_2\circ N_{E/F}$ | $I_\alpha(s_1-1/2,\delta(\chi_1)\otimes\chi_2^{-2})$ | $I_\alpha(s_1-1/2,\chi_1\circ det\otimes\chi_2^{-2})$ |
| $s_1=1$ and $\chi_1=1$ | $I_\alpha(s_2-3/2,\delta(\chi_2\circ N_{E/F})\otimes\chi_2^{-2})$ | $I_\alpha(s_2-3/2,(\chi_2\circ N_{E/F})\circ det\otimes\chi_2^{-2})$ |
| $3s_1-s_2=1$ and $\chi_1=\chi_2$ | $I_\beta(s_2-1/2,\chi_1\chi_2^{-1}\circ N_{E/F}\otimes\delta(\chi_2))$ | $I_\beta(s_2-1/2,\chi_1\chi_2^{-1}\circ N_{E/F}\otimes\chi_2\circ det)$ |
| $s_2=1$ and $\chi_2=1$ | $I_\beta(3s_1-3/2,\chi_1\chi_1^{-1}\circ N_{E/F}\otimes\delta(\chi_1))$ | $I_\beta(3s_1-3/2,\chi_1\chi_1^{-1}\circ N_{E/F}\otimes\chi_1\circ det)$ |

TABLE 3. Regular $\#R=1$

| Regular $\#R = 2$ | $I(s_1, s_2, \chi_1, \chi_2)$ | | | |
|---|---|---|---|---|
| | subrepresentation | | quotient | |
| | subrepresentation | quotient | subrepresentation | Langlands quotient |
| $(1, 2, \chi_1, 1; \chi_1\|_{F^\times} = 1)$ | $\pi(\chi_1) \simeq \pi(\chi_1^{-1})$ | $J_\beta(3/2, \chi_1 \otimes \delta(1))$ | $J_\beta(3/2, \chi_1^{-1} \otimes \delta(1))$ | $J_\alpha(1, I^\alpha(\chi_1, \chi_1^{-1}) \otimes 1)$ |
| $(1, 2, 1, \chi_2; \chi_2^2 = 1)$ | $\pi(\chi_2)$ | $J_\alpha(1/2, \delta(\chi_2 \circ N_{E/F}) \otimes 1)$ | $J_\beta(3/2, 1 \otimes \delta(\chi_2))$ | $J_\alpha(1, I^\alpha(1, \chi_2 \circ N_{E/F}) \otimes \chi_2)$ |
| $(1, 2, \chi_1, \chi_2; \chi_1^2 = 1, \chi_1 = \chi_2)$ | $\simeq (1, 2, 1, \chi_2; \chi_2^2 = 1)$ | | | |
| $(2, 3, 1, \chi_2; \chi_2 \circ N_{E/F} = 1)$ | $\pi(\chi_2) \simeq \pi(\chi_2^{-1})$ | $J_\alpha(3/2, \delta(1) \otimes \chi_2^{-1})$ | $J_\alpha(3/2, \delta(1) \otimes \chi_2)$ | $J_\beta(3, 1 \otimes I^\beta(\chi_2 \otimes \chi_2^{-1}))$ |
| $(3, 5, 1, 1)$ | $St_{G_2}$ | $J_\beta(9/2, 1 \otimes \delta(1))$ | $J_\alpha(5/2, \delta(1) \otimes 1)$ | $1_{G_2}$ |

Table 4. Regular $\#R = 2$

| Singular $1 \le \#R \le 2$ | $I(s_1, s_2, \chi_1, \chi_2)$ | |
|---|---|---|
| | subrepresentation | Langlands quotient |
| $(1, 2, 1, \chi_2; \chi_2 \circ N_{E/F} = 1)$ | $I_\alpha(1/2, \delta(1) \otimes \chi_2)$ | $J_\alpha(1, I^\alpha(1 \otimes 1) \otimes \chi_2^{-1})$ |
| $(2/3, 1, \chi_1, 1; \chi_1\|_{F^\times} = 1)$ | $I_\beta(1/2, \chi_1 \otimes \delta(1))$ | $J_\beta(1, \chi_1^{-1} \otimes I^\beta(1 \otimes 1))$ |
| $(2, 3, 1, 1)$ | $I_\alpha(3/2, \delta(1) \otimes 1)$ | $J_\beta(3, , 1 \otimes I^\beta(1 \otimes 1))$ |
| $(1, 3/2, 1, \chi_2; \chi_2^2 = 1)$ | $I_\alpha(0, \delta(\chi_2 \circ N_{E/F}) \otimes 1)$ | $J_\beta(3/2, 1 \otimes I^\beta(\chi_2 \otimes \chi_2))$ |
| $(1/2, 1, \chi_1, 1; \chi_1^2 = 1)$ | $I_\beta(0, 1 \otimes \delta(\chi_1))$ | $J_\alpha(1/2, I^\alpha(\chi_1 \otimes \chi_1) \otimes 1)$ |

Table 5. Singular length 2

*Remark* **2.** *If $E/F$ is a non-Galois cubic field extension, the previous Langlands classification almost holds. The only difference is that $N_{E/F}(E^\times) = F^\times$ (cf. Norm Limitation Theorem* [Mil97, Theorem 3.16]*). That is to say $(2, 3, 1, \chi_2; \chi_2 \circ N_{E/F} = 1)$ and $(1, 2, 1, \chi_2; \chi_2 \circ N_{E/F} = 1)$ will not appear in Tables 4 and 5 respectively.*

## 3. Unitary dual

In this section, we would like to sort out the unitary dual from our previous Langlands classification for $PGSO_8^E$. To do so, we first classify the Hermitian dual which states that

*For $\nu \in (\mathfrak{a}_M^*)_+ := \{x \in X^*(A_M) \otimes_{\mathbb{Z}} \mathbb{R} : (x, \alpha) > 0, \forall \alpha \in \Delta^G \backslash \Delta^M\}$ and $\sigma$ tempered, the Langlands quotient $J_P(\sigma \otimes \nu)$ is Hermitian if and only if there exists $w \in W(G, A_M)$ such that $\sigma \simeq w.\sigma$ and $-\nu = w.\nu$.*

Applying this criterion to our setting, we have, for irreducible tempered representations $\sigma$ of $GL_2$,

When $\chi_2$ is unitary and $s > 0$, the Langlands quotient $J_\alpha(s, \sigma \otimes \chi_2)$ is Hermitian if and only if $w_{3\alpha+2\beta}.(\sigma \otimes \chi_2) \simeq \sigma \otimes \chi_2$, i.e.

$$\sigma = \delta(\chi_1) \otimes 1|_{\chi_1^2 = 1} \text{ or } I^\alpha(\chi_1, \chi_1^{-1}\chi_2 \circ N_{E/F}) \otimes \chi_2|_{\chi_1^2 = 1 = \chi_2^2} \text{ or } I^\alpha(\chi_1, \chi_1^{-1}) \otimes 1;$$

When $\chi_1$ is unitary and $s > 0$, the Langlands quotient $J_\beta(s, \chi_1 \otimes \sigma)$ is Hermitian if and only if $w_{2\alpha+\beta}.(\chi_1 \otimes \sigma) \simeq \chi_1 \otimes \sigma$, i.e.

$$\sigma = 1 \otimes \delta(\chi_2)|_{\chi_2^2 = 1} \text{ or } \chi_1^{-1} \otimes I^\beta(\chi_2, \chi_2^{-1}\chi_1)|_{\chi_1^2 = 1 = \chi_2^2} \text{ or } 1 \otimes I^\beta(\chi_2, \chi_2^{-1});$$

When $\chi_1, \chi_2$ are unitary and $\frac{3}{2}s_2 < 3s_1 < 2s_2$, the Langlands quotient $J(s_1, s_2, \chi_1, \chi_2)$ is Hermitian if and only if $\chi_1^2 = 1 = \chi_2^2$.

For those Hermitian representations, we have the following associated reducibility conditions. As the discrete case has been discussed in Corollary 2.1, here we only consider the tempered non-discrete case.

**Lemma 3.1.** *Keep the notions as before. For unitary characters $\chi_1, \chi_2$ and $s > 0$, we have*

*(i) For $\chi_1^2 = 1$ and $\chi_2^2 = 1$, $I_\alpha(s, I^\alpha(\chi_1, \chi_1^{-1}\chi_2 \circ N_{E/F}) \otimes \chi_2)$ reduces if and only if*

$$s = 1/2, \chi_2 = 1 \quad \text{or} \quad s = 1, \chi_1|_{F^\times} = 1 \quad \text{or} \quad s = 1, \chi_1 = \chi_2.$$

*(ii) $I_\alpha(s, I^\alpha(\chi_1, \chi_1^{-1}) \otimes 1)$ reduces if and only if*

$$s = 1/2 \quad \text{or} \quad s = 1, \chi_1|_{F^\times} = 1.$$

*(iii)* *For* $\chi_1^2 = 1$ *and* $\chi_2^2 = 1$, $I_\beta(s, \chi_1^{-1} \otimes I^\beta(\chi_2, \chi_1\chi_2^{-1}))$ *reduces if and only if*

$$s = 3/2, \chi_1 = 1 \quad or \quad s = 1, \chi_1 = \chi_2 \text{ or } \chi_2 = 1 \quad or \quad s = 3, \chi_1 = \chi_2 \text{ or } \chi_2 = 1.$$

*(iv)* $I_\beta(s, 1 \otimes I^\beta(\chi_2, \chi_2^{-1}))$ $(\chi_2 \neq 1)$ *reduces if and only if*

$$s = \frac{3}{2} \quad or \quad s = 3, \chi_2 \circ N_{E/F} = 1.$$

In order to detect the unitarizability, we need to introduce another key input developed by Tadić and Speh, and summarized by Muić [Mui98, Lemma 5.1]. For an $F$-parabolic subgroup $P = MN$ of $G$, we denote by the $Unr(M)$ the group of unramified characters. For any irreducible representation $\sigma$ of $M$ and $\chi \in Unr(M)$, denote $I(\chi, \sigma) = Ind_P^G(\chi \otimes \sigma)$.

**Lemma 3.2.** *(*[Mui98, Lemma 5.1]*) Under the above assumptions, we have*

*(i)* *The set of those* $\chi \in Unr(M)$, *such that* $I(\chi, \sigma)$ *has a unitarizable irreducible subquotient, is compact.*

*(ii)* *Let* $S \subset Unr(M)$ *be a connected set. Suppose that for all* $\chi \in S$, *the representation* $I(\chi, \sigma)$ *is an irreducible unitarizable representation. Then for* $\chi \in \bar{S}$ *the closure of* $S$, *any irreducible subquotient of* $I(\chi, \sigma)$ *is unitarizable.*

*(iii)* *Suppose that* $\sigma$ *is Hermitian, and* $I(1, \sigma)$ *is irreducible and unitarizable. Then* $\sigma$ *is unitarizable.*

Before proceeding to sort out the whole unitary dual, we first verify some special cases as follows.

**Lemma 3.3.** *Suppose that* $\chi_1, \chi_2$ *are quadratic unitary characters and* $s > 0$. *Then*

$I_\alpha(s, \chi_1 \circ det \otimes 1)$ *is unitarizable (away from points of reducibility) if and only if* $s < 1/2$,

*and*

$I_\beta(s, 1 \otimes \chi_2 \circ det)$ *is unitarizable (away from points of reducibility) if and only if* $s < 3/2$.

*Proof.* This follows from the same argument as in [Mui98, Lemma 5.2]. □

Now let us turn to determine the unitary dual of $PGSO_8^E$. By Corollary 2.1 and Lemmas 3.1 and 3.2, we have

**Theorem 3.4.** *Keep the notation as before. For* $\chi_1, \chi_2$ *unitary characters of* $F^\times$ *and* $s > 0$, *we have*

*(i)* *For* $\chi_1^2 = 1$, $J_\alpha(s, \delta(\chi_1) \otimes 1)$ *is unitarizable if and only if* $s \le 1/2$.

*(ii)* *For* $\chi_2^2 = 1$, $J_\beta(s, 1 \otimes \delta(\chi_2))$ *is unitarizable if and only if* $s \le 3/2$.

*(iii)* $J_\alpha(s, I^\alpha(\chi_1, \chi_1^{-1}) \otimes 1)$ *is unitarizable if and only if* $s \le 1/2$, *or* $\chi_1|_{F^\times} = 1$ *and* $s = 1$.

*(iv)* *For* $\chi_1^2 = 1$ *and* $\chi_2^2 = 1$ $(\chi_2 \neq 1)$, $J_\alpha(s, I^\alpha(\chi_1, -) \otimes \chi_2)$ *is unitarizable if and only if* $\chi_1 = 1$ *and* $s \le 1$, *or* $\chi_1 = \chi_2$ *and* $s \le 1$.

*(v)* *For* $\chi_1^2 = 1$ *and* $\chi_2^2 = 1$ $(\chi_1 \neq 1)$, $J_\beta(s, \chi_1 \otimes I^\beta(\chi_2, -))$ *is unitarizable if and only if* $\chi_2 = 1$ *and* $s \le 1$, *or* $\chi_1 = \chi_2$ *and* $s \le 1$.

*(vi)* $J_\beta(s, 1 \otimes I^\beta(\chi_2, \chi_2^{-1}))$ *is unitarizable if and only if* $s \le \frac{3}{2}$, *or* $s = 3$ *provided* $\chi_2 \circ N_{E/F} = 1$ *and* $\chi_2 \neq 1$.

*Proof.* Following the standard procedure to construct families of positive definite Hermitian forms as in [Mui98, Theorem 5.1], we have

*Proof of (i)(ii).* It suffices to show $J_\alpha(5/2, \delta(\chi_1) \otimes 1)$ and $J_\beta(9/2, 1 \otimes \delta(\chi_2))$ are non-unitarizable which is well-known (see [BW13, Chapter XI Theorem 4.5]).

*Proof of (iii)(iv).* It suffices to show

$$I_\alpha(s, I^\alpha(\chi_1, \chi_1^{-1}) \otimes 1) \text{ is non-unitarizable for some } s \in (1/2, 1)$$

which follows from the fact that $I_\alpha(s, I^\alpha(\chi_1, \chi_1^{-1}) \otimes 1) = I_\beta(\chi_1 \otimes I^\beta(\nu_F^s \otimes \nu_F^{-s}))$ and Lemma 3.2 (iii), and

$$(\star) \qquad J_\alpha(1, I^\alpha(\chi_1, \chi_1^{-1}) \otimes 1) \text{ is unitarizable.}$$

If $\chi_1 = 1$, $(\star)$ follows from Lemma 3.3 and the fact that

$$I_\alpha(1/2, 1_{GL_2} \otimes 1) \twoheadrightarrow J_\alpha(1, I^\alpha(1, 1) \otimes 1),$$

if $\chi_1 \neq 1$, $(\star)$ will be proved later on.

*Proof of (v).* It suffices to show

$$I_\beta(s,\chi_1^{-1}\otimes I^\beta(1,1)) \text{ and } I_\beta(s,\chi_1^{-1}\otimes I^\beta(\chi_2,1))|_{\chi_1=\chi_2} \text{ are non-unitarizable for some } s\in(1,3)$$

which results from the fact that they are isomorphic to $I_\alpha(I^\alpha(\nu_E^{\frac{1}{3}s}\chi_1\otimes\nu_E^{-\frac{1}{3}s})\otimes\chi_1^{-1})$ and Lemma 3.2 (iii), and

$$J_\beta(2,3,\chi_1,1)\simeq I_\alpha(3/2,\chi_1\circ det\otimes 1) \text{ and } J_\beta(2,3,\chi_1,\chi_1)\simeq I_\alpha(3/2,\chi_1\circ det\otimes 1) \text{ are non-unitarizable}$$

which is considered in Lemma 3.3.

*Proof of (vi).* It suffices to show

$$I_\beta(\frac{2}{3}s,s,1,\chi_2)\simeq I_\alpha(I^\alpha(\nu_E^{s/3}\otimes\nu_E^{-s/3})\otimes\chi_2^{-1}) \text{ is non-unitarizable for some } s\in(3/2,3)$$

which is known by Lemma 3.2 (iii), and

$$I_\beta(\frac{2}{3}s,s,1,1)\simeq I_\alpha(I^\alpha(\nu_E^{\frac{1}{3}s}\otimes\nu_E^{-\frac{1}{3}s})\otimes 1) \text{ is unitarizable for } s\in(1,3/2)$$

which is also known by Lemma 3.2 (iii), and

$$J_\beta(3,1\otimes I^\beta(1\otimes 1))\simeq I_\alpha(3/2,1_{GL_2}\otimes 1) \text{ is non-unitarizable}$$

which is known by Lemma 3.3, and

$$J_\beta(3,1\otimes I^\beta(\chi_2,\chi_2^{-1})) \text{ is unitarizable provided } \chi_2\circ N_{E/F}=1 \text{ and } \chi_2\neq 1$$

which will be proved later on.

□

Before heading to the last case of unitarizable non-tempered Langlands quotients supported on the minimal parabolic subgroup, we recall the associated reducibility conditions as usual in the following.

**Lemma 3.5.** *For quadratic unitary characters $\chi_1,\chi_2$, and $(s_1,s_2)\in C^+$ the positive Weyl chamber, i.e. $\frac{1}{2}s_2<s_1<\frac{2}{3}s_2$. We know that $I(s_1,s_2,\chi_1,\chi_2)$ reduces if and only if $(s_1,s_2,\chi_1,\chi_2)$ is one of the following.*

$$(s_1,1,\chi_1,1;1/2<s_1<2/3),\quad (s_1,2s_1-1,\chi_1,1;s_1>2),\quad (1,s_2,1,\chi_2;3/2<s_2<2),$$
$$(\frac{2s_2-1}{3},s_2,1,\chi_2;s_2>2),\quad (s_1,3s_1-1,\chi_1,\chi_1;2/3<s_1<1),\quad (s_1,s_1+1,\chi_1,\chi_1;1<s_1<2).$$

**Theorem 3.6.** *Suppose that $\chi_1$ and $\chi_2$ are unitary characters, and $s_1$ and $s_2$ satisfy the condition that $\frac{1}{2}s_2<s_1<\frac{2}{3}s_2$, then $J(s_1,s_2,\chi_1,\chi_2)$ is unitarizable if and only if one of the following conditions holds.*

*(i)* $\chi_1=\chi_2=1$, *and* $s_2\leq 1$ *or* $3s_1-s_2\geq 1, s_2-s_1\leq 1$ *or* $s_1=3,s_2=5$.
*(ii)* $\chi_1=1,\chi_2$ *is of order 2, and* $s_1\leq 1$.
*(iii)* $\chi_2=1,\chi_1$ *is of order 2, and* $s_2\leq 1$.
*(iv)* $\chi_1=\chi_2$ *is of order 2, and* $3s_1-s_2\leq 1$.

*Proof.* By Lemma 3.5, we only have to discuss four cases as follows.

(i) $\underline{\chi_1=\chi_2=1}$: This results from the following analysis on the bounded domains partitioned by the reducibility lines case by case.

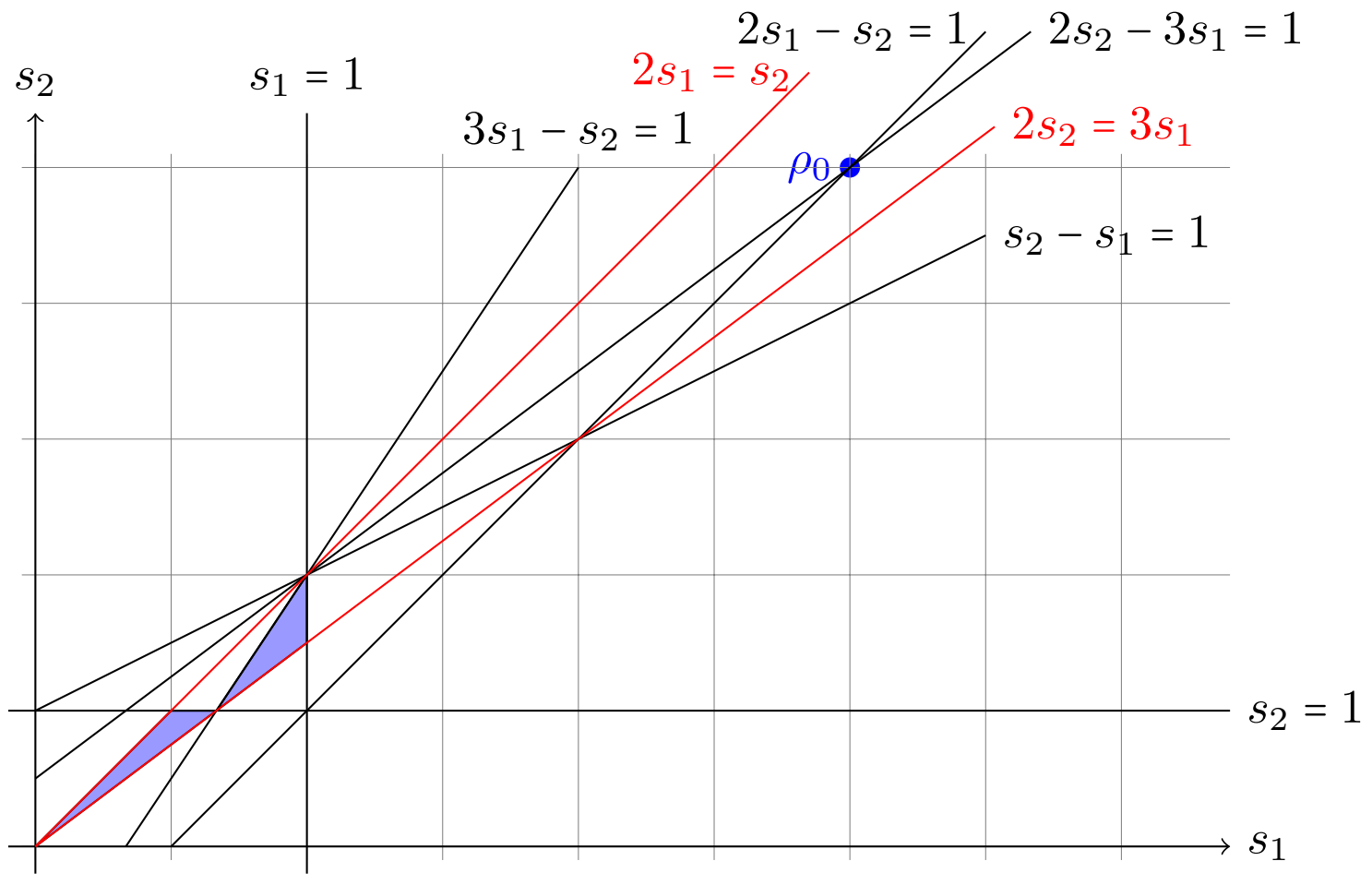

FIGURE 1. Spherical unitary dual with real infinitesimal character

(i1) $\underline{s_2 \leq 1}$: This is because $I(1 \otimes 1)$ is unitarizable.

(i2) $\underline{s_2 > 1, 3s_1 - s_2 < 1}$ : As $I(s_1, 2s_1, 1, 1) \simeq I_\alpha(s_1, I^\alpha(1 \otimes 1) \otimes 1)$ is non-unitarizable for $s_1 \in (1/2, 1)$ by Theorem 3.4 (iii).

(i3) $\underline{3s_1 - s_2 \geq 1, s_1 \leq 1}$ : It suffices to prove the unitarizability of one of the representations $J(s_1, s_2, 1, 1)$ under the condition $3s_1 - s_2 > 1, s_1 < 1$ by Lemma 3.2 (ii). The argument is the same as in the proof part (i3) of [Mui98, Theorem 5.2] by replacing $\beta$ by $\alpha$. But for completeness, we write down the argument as follows.
Consider

$$X_{t,s} := I_\alpha(t, I^\alpha(\nu_E^s \otimes \nu_E^{-s}) \otimes 1) = I(t + s, 2t, 1, 1).$$

The idea is to show the existence of an irreducible unitarizable domain $U$ of $X_{t,s}$ such that

$$(\star) \qquad \{(t + s, 2t)^w : w \in W,\ (t, s) \in U\} \cap \{(s_1, s_2) \in C^+ : 3s_1 - s_2 > 1,\ s_1 < 1\} \neq \varnothing.$$

We first classify the reducibility lines of $X_{t,s}$ as follows.

$$s = \frac{1}{2}, \quad t \pm 3s = 1, \quad t \pm s = 1, \quad t = \frac{1}{2}.$$

Then we sort out an irreducible unitarizable domain $U$ of $X_{t,s}$

$$U =: \{(t, s) : s \in (\frac{1}{3}, \frac{1}{2}),\ t \in (0, 1 - s)\}$$

It is quite easy to check that such a $U$ satisfying the above requirement $(\star)$.

(i4) $\underline{s_1 > 1, s_2 - s_1 < 1}$: As $I(s_1, \frac{3}{2}s_1, 1, 1) \simeq I_\beta(\frac{3}{2}s_1, 1 \otimes I^\beta(1 \otimes 1))$ is non-unitarizable for $1 < s_1 < 2$ by Theorem 3.4 (vi).

(i5) $\underline{s_2 - s_1 \geq 1, 2s_2 - 3s_1 < 1, 2s_1 - s_2 > 1}$ : On the boundary $s_2 - s_1 = 1$ with $1 < s_1 < 2$, we know the non-unitarizability of $J(s_1, s_1 + 1, 1, 1)$ which follows from the fact that

$$J(s_1, s_1 + 1, 1, 1) \simeq I_\alpha(s_1 - 1/2, 1_{GL_2} \otimes 1) \text{ is non-unitarizable for } 1 < s_1 < 2 \text{ by Lemma 3.3.}$$

(i6) $\underline{2s_2 - 3s_1 = 1, 1 < s_1 < 3}$ : As was known,

$$J(\frac{2s_2 - 1}{3}, s_2, 1, 1) \simeq I_\beta(s_2 - 1/2, 1 \otimes 1_{GL_2}) \text{ is non-unitarizable for } s_2 \in (2, 5) \text{ by Lemma 3.3.}$$

(i7) $\underline{2s_1 - s_2 = 1, 2 < s_1 < 3}$ : Similarly, this follows from the fact that

$$J(s_1, 2s_1 - 1, 1, 1) \simeq I_\alpha(s_1 - 1/2, 1_{GL_2} \otimes 1) \text{ is non-unitarizable for } s_1 \in (2, 3) \text{ by Lemma 3.3.}$$

(i8) $\underline{s_1 = 3, s_2 = 5}$ : $J(3, 5, 1, 1) \simeq 1_G$ is a unitarizable representation.

(ii) $\chi_1 = 1, \chi_2$ order 2: This follows from the fact that there is only one connected bounded domain determined by the reducibility lines.

(iii) $\chi_2 = 1, \chi_1$ order 2: This follows from the fact that

$$J(s_1, 2s_1 - 1, \chi_1, 1) \simeq I_\alpha(s_1 - 1/2, \chi_1 \circ det \otimes 1) \text{ is non-unitarizable for } s_1 > 2 \text{ by Lemma 3.3.}$$

(iv) $\chi_1 = \chi_2$ order 2: This follows from the fact that

$$J(s_1, s_1 + 1, \chi_1, \chi_1) \simeq I_\alpha(s_1 - 1/2, \chi_1 \circ det \otimes 1) \text{ is non-unitarizable for } s_1 > 1 \text{ by Lemma 3.3.}$$

□

**Unitary dual supported on $P_\gamma$.** Let $K$ be a non-archimedean field of characteristic zero, and denote by $W_K$ the associated Weil group of $K$. Let $\rho = \pi(\tau)$ be any supercuspidal representation of $GL_2(K)$, where

$$\tau : W_K \longrightarrow GL_2(\mathbb{C})$$

is an attached irreducible admissible homomorphism. Then $det\ \tau = \omega_\rho$ (the central character of $\rho$) via class field theory (see [Sha89, Section 1] for the details).

**Theorem 3.7.** *Suppose that $\rho$ is a unitary supercuspidal representation of $GL_2(K)$ for $K = F$ or $E$. We have*

*(i) The Langlands quotient $J_\alpha(s, \rho \otimes \chi_2)$ provided $\omega_\rho \chi_2 \circ N_{E/F} = 1$ is unitarizable if and only if $\rho \simeq \tilde{\rho}$ (the contragredient) and one of the following conditions holds:*
- *$\chi_2 = 1$ and $0 < s \le 1/2$.*
- *$0 < s \le 1$ and $\rho = Ind_{W_{E^c}}^{W_E}(\chi_0)$ provided $\chi_0|_S = 1$ and $\chi_2 \circ N_{S/F} = 1$, where $E^c/F$ is a Galois extension of degree 6 and $S \subset E^c$ is the unique quadratic extension over $F$.*

*(ii) The Langlands quotient $J_\beta(s, \chi_1 \otimes \rho)$ provided $\omega_\rho \chi_1 = 1$ is unitarizable if and only if $\rho \simeq \tilde{\rho}$ and one of the following conditions is satisfied:*
- *$\chi_1 = 1$ and $0 < s \le 1/2$.*
- *$Im(\tau) \simeq S_3$ (the symmetric group) given by the non-Galois extension $E$ over $F$, and $0 < s \le 1$.*

*If $I_\alpha(s_0, \rho \otimes \chi_2)$ (resp. $I_\beta(s_0, \chi_1 \otimes \rho)$), $s_0 > 0$, reduces, then it has a unique irreducible subrepresentation $\pi_\alpha(s_0, \rho \otimes \chi_1)$ (resp. $\pi_\beta(s_0, \chi_1 \otimes \rho)$). Those subrepresentation are square-integrable and different ($s_0$ is uniquely determined by the pair $(\rho, \chi_i)$). If $I_\alpha(0, \rho \otimes \chi_2)$ (or $I_\beta(0, \chi_1 \otimes \rho)$) reduces, then it is of length 2 and multiplicity 1.*

*Proof.* This follows from the $L$-factor computation in [Sha88, KK11] and the recent result of Henniart and Lomelí [HL17]. For $M_\alpha \simeq GL_2(E) \times F^\times / \Delta E^\times$, we have, using the standard notation as in [Sha90],

$$L(s, \rho \otimes \chi_2, r_2) = L_F(s, \chi_2^{-1})$$

and

$$L(s, \rho \otimes \chi_2, r_1) = L_F(s, \bigotimes -Ind_{W_E}^{W_F}(\tau) \otimes det^{-1}(\tau)) \ (\bigotimes -Ind \text{ twisted tensor induction [GS14, §6.1]}).$$

In view of those and the poles of twisted local triple product L-function which is proved in the Appendix (see also [Ike92, Theorem 2.6]), part (i) holds. As for $M_\beta \simeq E^\times \times GL_2(F)/\Delta F^\times$, we have

$$L(s, \chi_1 \otimes \rho, r_2) = L_E(s, \chi_1^{-1}), \quad L(s, \chi_1 \otimes \rho, r_3) = L_F(s, \tau)$$

and

$$L(s, \chi_1 \otimes \rho, r_1) = L_E(s, \chi_1 \cdot \rho_E)$$

where $\rho_E$ is the base change of $\rho$. In view of those, part (ii) holds. □

*Remark* **3.** *For the non-Galois cubic extension $E/F$ case, there is a new family of unitary representations concerning part (i) of Theorem 3.7 under the conditions that $0 < s \le 1$ and $\tau|_{W_{E^c}} = Ind_{W_L}^{W_{E^c}}(\chi_0)$ is irreducible, where $L/F$ is a Galois extension with $Gal(L/F) = D_{12}$ and $E^c/F$ is the Galois closure of $E/F$, such that*
- *$\chi_0|_{S^\times} \cdot \chi_2 \circ N_{S/F} = 1$, where $S \subset L$ is the degree 4 extension over $F$.*
- *$\omega_\rho \cdot \chi_2 \circ N_{E/F} = 1$.*
- *$\omega_\rho \circ N_{E^c/E} = \chi_0|_{(E^c)^\times} \cdot \omega_{L/E^c}$, where $\omega_{L/E^c}$ is the quadratic character associated to $L/E^c$.*

Note that J. Bernstein's unitarity conjecture says that the Aubert duality preserves unitarity. Back to our $PGSO_8^E$-setting, based on our computation, we have

**Corollary 3.8.** *Keep the notation as before. The unitary dual is preserved under the Aubert duality.*

Note also that L. Clozel's finiteness conjecture (see [Clo85] for the details) says that the set of exponents of discrete series is finite. Put in our setting, we have

**Corollary 3.9.** *Keep the notions as before. Clozel's finiteness conjecture of special exponents holds for $PGSO_8^E$.*

**Unitarizability of $J_\alpha(1, I^\alpha(\chi_1, \chi_1^{-1}) \otimes 1)$ and $J_\beta(3, 1 \otimes I^\beta(\chi_2, \chi_2^{-1}))$.** In what follows, we prove that $J_\alpha(1, I^\alpha(\chi_1, \chi_1^{-1}) \otimes 1)$ (resp. $J_\beta(3, 1 \otimes I^\beta(\chi_2, \chi_2^{-1})))$, where $\chi_1|_{F^\times} = 1$ and $\chi_1 \neq 1$ (resp. $\chi_2 \neq 1$ and $\chi_2 \circ N_{E/F} = 1$), is a unitary representation. Then it is an isolated point in the unitary dual of $PGSO_8^E$ by [Tad88, Theorem 2.2].

The main idea is to show that they appear as components of some specific residual spectrum of $G$ as in [Mui98, Žam97, Kim96]. Let us start with some notation. For a global field $\dot{K}$, let $\mathbb{A}_{\dot{K}}$ be the ring of Adeles of $\dot{K}$. As in the local field case, given $\dot{E}$ a cubic field extension of a global field $\dot{F}$, we have an associated quasi-split adjoint group $G = PGSO_8^{\dot{E}}$ of type $D_4$. For grössencharacters $\mu_1$ and $\mu_2$ of $\dot{E}$ and $\dot{F}$ respectively, we define a unitary character $\chi = (\mu_1, \mu_2)$ of $T(\mathbb{A}_{\dot{F}})$ by $\chi(t(a,b)) = \mu_1(a)\mu_2(b)$. We take the coordinates in $\mathfrak{a}_{\mathbb{C}}^* = X^*(T) \otimes \mathbb{C}$ with respect to the basis $\alpha$, $\beta$; the ordered pair $(s_1, s_2) \in \mathbb{C}^\times$ corresponds to the character $\lambda = 3s_1\alpha + s_2\beta$. For $\lambda$ and $\chi$ as above, let $I_B(\lambda, \chi) = I_{T(\mathbb{A}_{\dot{F}})}^{G(\mathbb{A}_{\dot{F}})}(\lambda, \chi)$ be the space for the standard normalized induction (sometimes write as $I_B(\nu_{\dot{E}}^{s_1}\mu_1 \otimes \nu_{\dot{F}}^{s_2}\mu_2)$). Finally, let $\rho_B$ be the half sum of positive roots, i.e. $\rho_B = 5\alpha + 3\beta$, and $C^+$ be the positive Weyl chamber in $\mathfrak{a}_{\mathbb{C}}^*$:

$$C^+ = \{s_1\alpha + s_2\beta : \frac{3}{2}Re(s_2) < Re(s_1) < 2Re(s_2)\}.$$

Following the standard procedure of investigating $L_d^2(B)$,

- (Eisenstein series) For $f \in I_B(\lambda, \chi)$, one forms Eisenstein series

$$E(g, f, \lambda) = \sum_{\gamma \in B(\dot{F}) \backslash G(\dot{F})} f(\gamma g)$$

  which converges absolutely for $Re\lambda \in C^+ + \rho_B$ and extends to a meromorphic function of $\lambda$. It is an automorphic form and its singularities coincide with those of its constant term along $B$, i.e.

(C) $$E_0(g, f, \lambda) = \sum_{w \in W} (M(w, \lambda, \chi) f)(g)$$

  where $M(w, \lambda, \chi) = \otimes_\nu M(w, \lambda, \chi_\nu)$ are the so-called non-normalized intertwining operators from $I_B(\lambda, \chi)$ to $I_B(w\lambda, w\chi)$.
- (Normalization) Let $\psi = \otimes_\nu \psi_\nu$ be a fixed non-trivial additive character of $\dot{F} \backslash \mathbb{A}_{\dot{F}}$. The standard normalization of the intertwining operators $M(w, \lambda, \chi)$ for all $\nu$ by factors (assume $\dot{E}/\dot{F}$ is Galois for simplicity),

$$r(w, \lambda, \chi_\nu) = \prod_{\substack{\{\gamma > 0, w.\gamma < 0\} \\ \gamma \ long}} \frac{L(\langle \lambda, \gamma^\vee \rangle, \chi_\nu \circ \gamma^\vee)}{L(\langle \lambda, \gamma^\vee \rangle + 1, \chi_\nu \circ \gamma^\vee) \epsilon(\langle \lambda, \gamma^\vee \rangle, \chi_\nu \circ \gamma^\vee, \psi_\nu)}$$
$$\times \prod_{\substack{\{\gamma > 0, w.\gamma < 0\} \\ \gamma \ short}} \frac{L(\langle \lambda, \gamma^\vee \rangle / 3, \chi_\nu \circ \gamma^\vee)}{L(\langle \lambda, \gamma^\vee \rangle / 3 + 1, \chi_\nu \circ \gamma^\vee) \epsilon(\langle \lambda, \gamma^\vee \rangle / 3, \chi_\nu \circ \gamma^\vee, \psi_\nu)}$$

  are as follows:

$$N(w, \lambda, \chi_\nu) = r(w, \lambda, \chi_\nu)^{-1} M(w, \lambda, \chi_\nu) \text{ which are multiplicative.}$$

  Let $N(w, \lambda, \chi) = \otimes_\nu N(w, \lambda, \chi_\nu)$, it is well-known that

$$M(w, \lambda, \chi_\nu) \prod_{\substack{\{\gamma > 0, w.\gamma < 0\} \\ \gamma \ long}} L(\langle \lambda, \gamma^\vee \rangle, \chi_\nu \circ \gamma^\vee)^{-1} \prod_{\substack{\{\gamma > 0, w.\gamma < 0\} \\ \gamma \ short}} L(\langle \lambda, \gamma^\vee \rangle / 3, \chi_\nu \circ \gamma^\vee)^{-1}$$

is holomorphic for all $\nu$.

- (Singularities) The possible singularities of $M(w,\lambda,\chi)$ are rank 1 reducibility points as in Table 1, zeros of the denominator of $r(w,\lambda,\chi)=\prod_\nu r(w,\lambda,\chi_\nu)$ and poles of $N(w,\lambda,\chi)$. It is easy to see that only the point $3\alpha+\beta$ could provide a pole of $N(w,\lambda,\chi)$ as in [Žam97].
- (Langlands square-integrable criterion) $Res_{\lambda_0}Res_{\langle\lambda,\gamma^\vee\rangle=1}E(g,f,\lambda)$ is square-integrable if and only if
$$Re(w\lambda_0)\in\{-u\alpha-v\beta:u,v>0\},\ \text{for all } w\in W_0$$
where $W_0\subset W$ consists of those elements that give non-zero residue on the right-hand side of (C) which is non canceled by residue of any other term.

$J_\beta(3,1\otimes I^\beta(\chi_2,\chi_2^{-1}))$ unitary: This results from the following two lemmas as in [Mui98, Theorem 6.2].

**Lemma 3.10.** *Let $\chi$ be a grössencharacter of $\dot{F}$ of order 3. Then the representation*
$$J_\beta(3,1\otimes I^\beta(\chi,\chi^{-1}))=J_\beta(3,1\otimes I^\beta(\chi^{-1},\chi))=\bigotimes_\nu J_\beta(3,1\otimes I^\beta(\chi_\nu,\chi_\nu^{-1}))$$
*occurs in the residual spectrum of $G$.*

*Proof.* This is to take residue at $\Lambda=6\alpha+3\beta$. It is easy to see that the point $\Lambda=6\alpha+3\beta$ only gives rise to simple poles arising from $r(w,\Lambda,\chi)$. So $W_0\subset W_{1,5}:=\{w\in W: w.\alpha<0, w.(\alpha+\beta)<0\}=\{w_{2\alpha+\beta},w_\beta w_{2\alpha+\beta}\}$. The same argument as in [Žam97, Case a) Residue at $\Lambda=2\alpha+\beta$], we know that $W_0=W_{1,5}$ and the residue of the constant term (C) produces $J_\beta(3,1\otimes I^\beta(\chi,\chi^{-1}))=\otimes_\nu J_\beta(3,1\otimes I^\beta(\chi_\nu,\chi_\nu^{-1}))$, whence the lemma holds. □

**Lemma 3.11.** ([AT67, Theorem 5]) *Let $\dot{K}$ be a global field, and $S$ be a finite set of places of $\dot{K}$. For $\nu\in S$, let $\chi_\nu$ be a character of $\dot{K}_\nu^\times$ of order dividing $n\in\mathbb{N}$. Then there exists a character $\mu$ of $\dot{K}\backslash\mathbb{A}_{\dot{K}}^\times$ of order dividing $2n$, such that $\mu_\nu=\chi_\nu$ for $\nu\in S$.*

$J_\alpha(1,I^\alpha(\chi_1,\chi_1^{-1})\otimes1)$ unitary: This results from the same argument as above.

**Lemma 3.12.** *Let $\chi$ be a grössencharacter of $\dot{E}$ such that $\chi|_{\mathbb{A}_{\dot{F}}^\times}=1$. Then the representation*
$$J_\alpha(1,I^\alpha(\chi,\chi^{-1})\otimes1)=\bigotimes_\nu J_\alpha(1,I^\alpha(\chi_\nu,\chi_\nu^{-1})\otimes1)$$
*occurs in the residual spectrum of $G$.*

*Proof.* This is about taking residue at $3\alpha+2\beta$. It is easy to see that the point $\Lambda=3\alpha+2\beta$ only gives rise to simple poles arising from $r(w,\Lambda,\chi)$. So $W_0\subset W_{2,6}:=\{w\in W: w.\beta<0, w.(3\alpha+\beta)<0\}=\{w_{3\alpha+2\beta},w_\alpha w_{3\alpha+2\beta}\}$. The same argument as in [Žam97, Case a) Residue at $\Lambda=2\alpha+\beta$], we know that $W_0=W_{2,6}$ and the residue of the constant term (C) produces $J_\alpha(1,I^\alpha(\chi,\chi^{-1})\otimes1)=\otimes_\nu J_\alpha(1,I^\alpha(\chi_\nu,\chi_\nu^{-1})\otimes1)$, whence the lemma holds. □

**Lemma 3.13.** *Let $\dot{E}$ be a cubic extension of a global field $\dot{F}$, and $S$ be a finite set of places of $\dot{F}$. For $\nu\in S$, let $\chi_\nu$ be a character of $\dot{E}_\nu$ such that $\chi_\nu|_{\dot{F}_\omega}=1$. Then there exists a grössencharacter character $\mu$ of $\dot{E}$, such that $\mu|_{\mathbb{A}_{\dot{F}}}=1$ and $\mu_\nu=\chi_\nu$ for $\nu\in S$.*

*Proof.* This follows from the fact:
$$\prod_{\nu\in S}\dot{E}_\nu^\times/\dot{F}_\nu^\times\longrightarrow\mathbb{A}_{\dot{E}}^\times/\dot{E}^\times\mathbb{A}_{\dot{F}}^\times \text{ is injective.}$$
□

## APPENDIX: POLES OF LOCAL TRIPLE PRODUCT L-FUNCTIONS

In this appendix, we will determine the poles of local triple product $L$-functions, which turns out to be the same as in the global case (treated by Ikeda in [Ike92]), but the proof is of course completely different, since the local proof proceeds on the Galois side based on the recent work of Henniart and Lomelí [HL17].

Let us first consider the case when $E = F\times F\times F$. Hence, let $\phi_1,\ \phi_2,\ \phi_3 :\ W_F\to GL_2(\mathbb{C})$ be three irreducible representations (corresponding to supercuspidal representations of $GL_2(F)$). We are interested in determining if $(\phi_1\otimes\phi_2\otimes\phi_3)^{W_F}\neq 0$. Equivalently, whether $\phi_1\otimes\phi_2$ can contain an irreducible 2-dimensional summand.

Suppose that $\phi_1\otimes\phi_2 = \rho_1\oplus\rho_2$ with $dim(\rho_i) = 2$.

Claim. $\phi_1$ and $\phi_2$ must have the form $\phi_i = Ind_{W_K}^{W_F}(\chi_i)$ for some quadratic field extension $K/F$ (independent of $i$), i.e. $\phi_1$ and $\phi_2$ are dihedral w.r.t. $K/F$.

Before justifying the claim, we first recall the following possibilities for $\phi := \phi_i$.

(a) $\phi$ is not dihedral

$$\Leftrightarrow \phi\bigotimes\chi\neq\phi \text{ for any quadratic character } \chi\neq 1$$
$$\Leftrightarrow \phi|_{W_K} \text{ is irreducible for any quadratic extension } K/F$$
$$\Leftrightarrow Sym^2\phi = \wedge^2\phi\bigotimes Ad(\phi) \text{ is irreducible.}$$

(b) $\phi$ is dihedral w.r.t. a unique quadratic extension $K/F$.

$$\Leftrightarrow \phi\bigotimes\omega_{K/F} = \phi,\ \textit{but}\ \phi\bigotimes\chi\neq\phi \text{ for any quadratic character } \chi\neq\omega_{K/F}\ \textit{or}\ 1.$$
$$\Leftrightarrow Ad(\phi)\ \textit{contains}\ \omega_{K/F}, \text{ but not other quadratic characters.}$$

In this case, we may write $\phi = Ind_{W_K}^{W_F}(\chi)$ for some character $\chi$ of $W_K$.

(c) $\phi$ is dihedral w.r.t. three quadratic extensions $K_i$ of $F$, $i = 1,2,3$.

$$\Leftrightarrow Ad(\phi) \text{ is the sum of three quadratic characters } \chi_1,\chi_2,\chi_3, \text{ such that } \chi_1\chi_2\chi_3 = 1.$$

In this case, we may write $\phi = Ind_{W_{K_i}}^{W_F}(\chi_i)$ for each $i$.

Now to justify the claim, we consider $\wedge^2$ on both sides of the equation $\phi_1\otimes\phi_2 = \rho_1\oplus\rho_2$. This gives:

$$(\star\star)\qquad (\wedge^2\phi_1\bigotimes Sym^2\phi_2)\bigoplus(Sym^2\phi_1\bigotimes\wedge^2\phi_2) = \wedge^2\rho_1\bigoplus\wedge^2\rho_2\bigoplus\rho_1\bigotimes\rho_2.$$

We now argue:

- At least one of $\phi_1,\phi_2$ is dihedral. If not, then $LHS$ of $(\star\star)$ is the sum of two 3-dimensional irreducible summands, whereas $RHS$ is not.
- If $\phi_1$ is dihedral, say $\phi_1 = Ind_{W_K}^{W_F}(\chi)$, but $\phi_2$ is not dihedral. Then
$$\phi_1\bigotimes\phi_2 = Ind_{W_K}^{W_F}(\chi\cdot\phi_2|_{W_K}).$$
Since $\phi_2|_{W_K}$ is irreducible, $\phi_1\otimes\phi_2$ is either irreducible, or a sum $\rho_1\oplus\rho_2 = \rho\oplus\rho\cdot\omega_{K/F}$. Looking at $(\star\star)$, one sees that $LHS$ contains either one or three distinct 1-dimensional characters, whereas $RHS$ contains $\wedge^2\rho_1 = \wedge^2\rho_2$ with multiplicity $\geq 2$.
- Thus both $\phi_1$ and $\phi_2$ are dihedral. If they are not dihedral w.r.t. the same $K$, then $\phi_1,\phi_2$ are as in case (b) above. Let $\phi_i = Ind_{W_{K_i}}^{W_F}(\chi_i)$. Then $\phi_1\otimes\phi_2 = Ind_{W_{K_1}}^{W_F}(\chi_1\cdot\phi_2|_{W_{K_1}})$ is either irreducible or the sum $\rho\oplus\rho\cdot\omega_{K/F}$. Looking at $(\star\star)$, we see that $LHS$ contains two distinct 1-dimensional characters, whereas $RHS$ contains $\wedge^2\rho_1 = \wedge^2\rho_2$ with multiplicity $\geq 2$.
- We have thus shown that there exists a quadratic extension $K/F$ such that $\phi_i = Ind_{W_K}^{W_F}(\chi_i)$. Then
$$\phi_1\bigotimes\phi_2 = Ind_{W_K}^{W_F}(\chi_1\chi_2)\bigoplus Ind_{W_K}^{W_F}(\chi_1\chi_2^\tau),$$
where $Gal(K/F) = \langle\tau\rangle$. Hence if $\tilde{\phi}_2$ (the contragredient) is a summand of $\phi_1\otimes\phi_2$, then $\tilde{\phi}_3$ is one of the two summands above, i.e. $\phi_3 = Ind_{W_K}^{W_F}(\chi_1\chi_2)^{-1}$ (replacing $\chi_2$ by $\chi_2^\tau$ if necessary).

We have shown:

**Proposition 3.14.** *Let $\phi_1,\ \phi_2,\ \phi_3 :\ W_F\to GL_2(\mathbb{C})$ be irreducible. Then $(\phi_1\otimes\phi_2\otimes\phi_3)^{W_F}\neq 0$ $\Leftrightarrow \exists$ quadratic extension $K/F$ s.t $\phi_i = Ind_{W_K}^{W_F}(\chi_i)$, with $\chi_i^\tau\neq\chi_i$ and $\chi_1\chi_2\chi_3 = 1$, in which case, the quadratic extension $K/F$ is uniquely determined by $\phi_1,\phi_2,\phi_3$ via:*

$$det(\phi_1)\cdot det(\phi_2)\cdot det(\phi_3) = \omega_{K/F}.$$

*Proof.* We have already shown the $(\Leftrightarrow)$. It remains to prove the last assertion.

With $\phi_i = Ind_{W_K}^{W_F}(\chi_i)$, $\chi_1\chi_2\chi_3 = 1$, one has $det(\phi_i) = \chi_i|_{F^\times} \cdot \omega_{K/F}$. So

$$det(\phi_1) \cdot det(\phi_2) \cdot det(\phi_3) = \chi_1\chi_2\chi_3|_{F^\times}\omega_{K/F}^3 = \omega_{K/F}.$$

□

*Remark* **4.** *As a consequence, we see that one cannot have $\phi_1$ and $\phi_2$ to be both dihedral w.r.t. the same three quadratic extensions $K_1, K_2, K_3$. For this will contradict the uniqueness part of the proposition.*

Now we consider the main case of interest where $E/F$ is a cubic field extension.

$E/F$ **Galois**. : We first consider the case that $E/F$ is a Galois extension. Suppose $Gal(E/F) = \langle\sigma\rangle$ and let $\tilde{\sigma} \in W_F$ be an element which projects to $\sigma$ under $W_F \twoheadrightarrow Gal(E/F)$. Let $\phi : W_E \to GL_2(\mathbb{C})$ be an irreducible representation and set

$$\rho = \bigotimes -Ind_{W_E}^{W_F}(\phi)$$

to be the tensor induction of $\phi$ from $W_E$ to $W_F$ (see [HL17, §2.1] for the notion of tensor induction), so that $dim(\rho) = 8$. We are interested in determining when $\rho^{W_F} \neq 0$.

Now $\rho^{W_F} \neq 0 \Rightarrow \rho^{W_E} \neq 0$. Since $\rho|_{W_E} = \phi \otimes \phi^\sigma \otimes \phi^{\sigma^2}$, our proposition shows that there exists a unique quadratic extension $L/E$ such that

$$\phi = Ind_{W_L}^{W_E}(\chi), \quad \phi^\sigma = Ind_{W_L}^{W_E}(\chi'), \quad \phi^{\sigma^2} = Ind_{W_L}^{W_E}(\chi'') \text{ with } \chi\chi'\chi'' = 1.$$

Claim: $L/F$ is a Galois extension.

*Proof.* It suffices to show that $\tilde{\sigma}(L) = L$. If not, then $L, \tilde{\sigma}(L), \tilde{\sigma}^2(L)$ are three distinct quadratic extensions of $E$. Moreover

$$\phi^\sigma = Ind_{W_L}^{W_E}(\chi') \Rightarrow \phi = Ind_{W_{\tilde{\sigma}^2(L)}}^{W_E}(\tilde{\sigma}^2(\chi'))$$
$$\phi^{\sigma^2} = Ind_{W_L}^{W_E}(\chi'') \Rightarrow \phi = Ind_{W_{\tilde{\sigma}(L)}}^{W_E}(\tilde{\sigma}(\chi'')).$$

So $\phi$ is dihedral w.r.t. $L$, $\tilde{\sigma}(L)$ *and* $\tilde{\sigma}^2(L)$. A similar argument shows the same for $\phi^\sigma$ and $\phi^{\sigma^2}$. This contradicts our earlier proposition, or rather the remark following it. So we must have $\tilde{\sigma}(L) = L$. □

As a consequence of the claim, $Gal(L/F) = \langle c\rangle$ is a cyclic group of order 6, and we have:

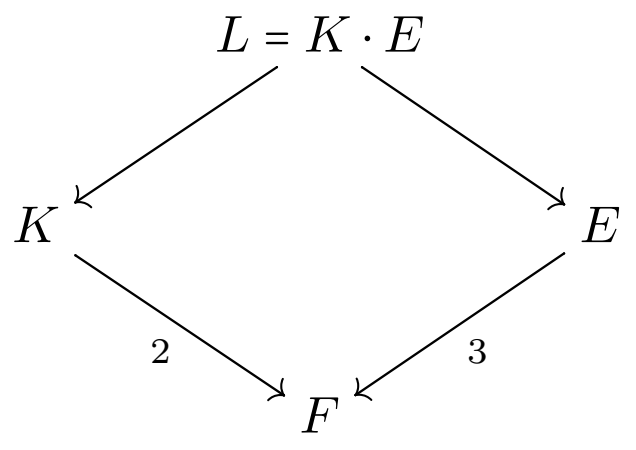

with $Gal(L/K) = \langle\tilde{\sigma}|_L\rangle = \langle c^2\rangle$ and $Gal(L/E) = \langle\tau\rangle = \langle c^3\rangle$.

Now a short computation shows (see [Ike92, Theorem 2.6]):

**Lemma 3.15.**

$$\rho := \bigotimes -Ind_{W_E}^{W_F}(Ind_{W_L}^{W_E}(\chi))$$
$$\simeq Ind_{W_K}^{W_F}(\chi|_{K^\times}) \bigoplus Ind_{W_L}^{W_F}(\chi^\tau \cdot \chi^{\tilde{\sigma}} \cdot \chi^{\tilde{\sigma}^2})$$

*where we have regarded $\chi$ as a character of $L^\times$.*

For $\rho^{W_F} \neq 0$, we need

$$\text{either} \quad \chi|_{K^\times} = 1 \qquad \text{or} \quad \chi^\tau\chi^{\tilde{\sigma}}\chi^{\tilde{\sigma}^2} = 1.$$

Let us show that the latter case is not possible. Indeed, if

$$1 = \chi^\tau\chi^{\tilde{\sigma}}\chi^{\tilde{\sigma}^2} = \chi^{c^3}\chi^{c^2}\chi^{c^4},$$

then applying $c$ gives:

$$1 = \chi^{c^4}\chi^{c^3}\chi^{c^5}.$$

Comparing the two equations gives:

$$\chi^{c^5} = \chi^{c^2},\ i.e.\ \chi^{c^3} = \chi,\ i.e.\ \chi^{\tau} = \chi.$$

But $\chi^{\tau} \neq \chi$ since $\phi$ is irreducible.

Hence we have shown:

**Theorem 3.16** ($E/F$ Galois)**.** *Let $\phi:\ W_E \to GL_2(\mathbb{C})$ be irreducible. Then $\rho := \otimes\text{-}Ind_{W_E}^{W_F}(\phi)$ contains the trivial character*

*$\Leftrightarrow \exists$ a quadratic extension $K/F$ and a character $\chi$ of $L^\times = (K{\cdot}E)^\times$ s.t $\phi = Ind_{W_L}^{W_E}(\chi)$ and $\chi|_{K^\times} = 1$, in which case,*

$$\rho \simeq Ind_{W_K}^{W_F}(1) \bigoplus Ind_{W_L}^{W_F}(\chi^{\tau}\chi^{-1})$$

*and $K$ is uniquely determined by*

$$\omega_{K/F} = \omega_{L/E}|_{F^\times} = det(\phi)|_{F^\times}.$$

<u>$E/F$ **non-Galois**</u>. Now we turn to the non-Galois case. Let $E^c/F$ be the Galois closure of $E/F$ with

$$Gal(E^c/F) = S_3 := \langle \tau, \sigma | \tau^2 = \sigma^3 = 1, \tau\sigma\tau = \sigma^{-1} \rangle.$$

We have two cases:

(i) <u>$\phi|_{W_{E^c}}$ reducible</u>: Similar argument as above shows that this gives rise to the same condition as in Theorem 3.16 for $\rho^{W_F} \neq 0$.

(ii) <u>$\phi|_{W_{E^c}}$ irreducible</u>: Similar argument as in the Galois case, $\rho^{W_{E^c}} \neq 0$ implies that there exists a unique quadratic extension $L/E^c$ such that

$$\phi|_{W_{E^c}} = Ind_{W_L}^{W_{E^c}}(\chi),\quad \phi^{\sigma}|_{W_{E^c}} = Ind_{W_L}^{W_{E^c}}(\chi'),\quad \phi^{\sigma^2}|_{W_{E^c}} = Ind_{W_L}^{W_{E^c}}(\chi'')\ with\ \chi\chi'\chi'' = 1.$$

Suppose $Gal(L/E^c) = \langle \tau' \rangle$. As $\phi|_{W_{E^c}}$ is irreducible, so

$$Ind_{W_L}^{W_{E^c}}(\chi) \simeq Ind_{W_{L^\tau}}^{W_{E^c}}(\chi^{\tau}),$$

which in turn implies that

$$L^{\tau} = L, \quad and \quad \chi = \chi^{\tau}\ or\ \chi^{\tau'} = \chi^{\tau}.$$

This is to say $L/E$ is a Galois extension and

$$Gal(L/E) \simeq \mathbb{Z}/2\mathbb{Z} \times \mathbb{Z}/2\mathbb{Z}.$$

Further applying Proposition 3.14, we know that $L^{\sigma} = L$, which in turn implies that $L/F$ is a Galois extension with

$$Gal(L/F) \simeq D_{12} := \langle \tau, \sigma_0 | \tau^2 = \sigma_0^6 = 1, \tau\sigma_0\tau = \sigma_0^{-1} \rangle,$$

and we have:

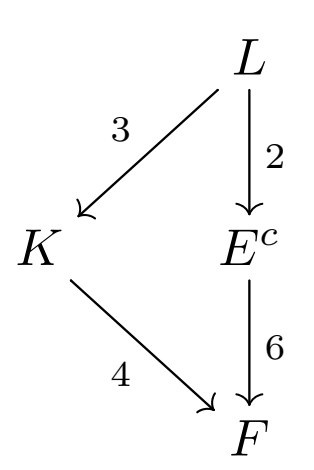

with $Gal(L/K) = \langle \sigma_0^2 \rangle$ and $Gal(K/F) = \langle \tau, \sigma_0^3 \rangle$.

Now a short computation shows that

$$\rho|_{W_K} = \chi\chi^{\sigma_0^2}\chi^{\sigma_0^4} + \chi^{\sigma_0}\chi^{\sigma_0^3}\chi^{\sigma_0^5} + Ind_{W_L}^{W_K}(\chi\chi^{\sigma_0}\chi^{\sigma_0^2}) + Ind_{W_L}^{W_K}(\chi\chi^{\sigma_0}\chi^{\sigma_0^5}).$$

So $\rho^{W_K} \neq 0$ implies that, applying the same argument as in Lemma 3.15,

$$\chi\chi^{\sigma_0^2}\chi^{\sigma_0^4}|_{W_K} = 1,\ i.e.\ \chi|_{K^\times} = 1.$$

On the other hand, given $\chi|_{K^\times} = 1$, an easy calculation shows that $\rho^{W_F} \neq 0$. Thus we obtain

**Theorem 3.17** ($E/F$ non-Galois)**.** *Let $\phi : W_E \to GL_2(\mathbb{C})$ be irreducible. Denote by $E^c/F$ the Galois closure of $E/F$. Then $\rho := \otimes - Ind_{W_E}^{W_F}(\phi)$ contains the trivial character if and only if one of the following conditions holds:*

*(i) There exists a character $\chi$ of $(E^c)^\times$, such that $\phi = Ind_{W_{E^c}}^{W_E}(\chi)$ and $\chi|_{K^\times} = 1$, here $K/F$ is the unique intermediate quadratic extension, in which case,*

$$\rho \simeq Ind_{W_K}^{W_F}(1) \bigoplus Ind_{W_L}^{W_F}(\chi^\tau \chi^{-1}).$$

*(ii) There exits a quadratic extension $L/E^c$ and a character $\chi$ of $L^\times$, such that $Gal(L/F) = D_{12}$, $\phi|_{W_{E^c}} = Ind_{W_L}^{W_{E^c}}(\chi)$ is irreducible and $\chi|_{K^\times} = 1$. Here $K/F$ is the unique quartic intermediate extension.*

*Remark* **5.** *As pointed out by Professor Ikeda, part (ii) of Theorem 3.17 is indeed a dihedral case given as follows.*

*As $Gal(L/E) \simeq \mathbb{Z}/2\mathbb{Z} \times \mathbb{Z}/2\mathbb{Z}$, we have the following diagram:*

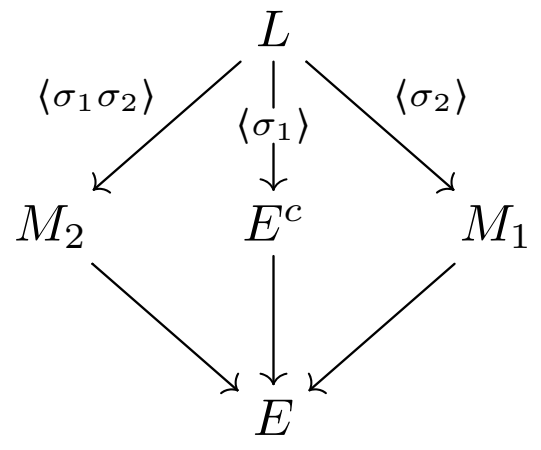

*The point is to show that*

$$\text{one of } Ind_{W_L}^{W_{M_i}}(\chi),\ i = 1, 2,\ \text{is irreducible.}$$

*Otherwise,*

(A) $$Ind_{W_L}^{W_{M_i}}(\chi) \text{ is irreducible for } i = 1, 2.$$

*That is to say*

$$\chi \not\simeq \chi^{\sigma_2} \text{ and } \chi \not\simeq \chi^{\sigma_1\sigma_2}.$$

*Note that*

(B) $$\phi|_{W_L} = \chi + \chi^{\sigma_1} \text{ with } \chi \not\simeq \chi^{\sigma_1} \ (\text{as } \phi|_{W_{E^c}} = Ind_{W_L}^{W_{E^c}}(\chi) \text{ is irreducible}).$$

*Therefore*

$$\begin{aligned} 0 \neq Hom_{W_L}(\chi, \phi) &= Hom_{W_E}(Ind_{W_L}^{W_E}(\chi), \phi) \\ &= Hom_{W_{M_i}}(Ind_{W_L}^{W_{M_i}}(\chi), \phi). \end{aligned}$$

*Thus (A) implies*

$$\phi|_{W_{M_i}} = Ind_{W_L}^{W_{M_i}}(\chi), \ for \ i = 1, 2.$$

*Which in turn says that*

$$\phi|_{W_L} = \chi + \chi^{\sigma_2} = \chi + \chi^{\sigma_1\sigma_2} \overset{(B)}{=} \chi + \chi^{\sigma_1}.$$

*Contradiction.*

Department of Mathematics, National University of Singapore, 10 Lower Kent Ridge Road Singapore 119076

*E-mail address*: cluo@u.nus.edu